\documentclass[preprint,11pt]{erkarticle}  

\usepackage{general,resistance}
\usepackage{pxfonts} %ifthen,leading}

\usepackage[bookmarks, colorlinks=true, pdfstartview=FitV, linkcolor=blue, citecolor=blue, urlcolor=blue]{hyperref}
\usepackage[left=1.6in,right=1.6in,top=1.7in,bottom=1.40in]{geometry}%showframe

\newcommand{\SubEn}{\negsp[7]\raisebox{-10pt}{\scriptsize \energy}}
\renewcommand{\linenopax}{} 

\numberwithin{equation}{section} \numberwithin{theorem}{section}

\begin{document}
%  \leading{13pt} %12pt
  %\pagewiselinenumbers

\begin{frontmatter}

\title{Self-adjoint extensions of network Laplacians and applications to resistance metrics}

%\thanks{The work of PETJ was partially supported by NSF grant DMS-0457581. The work of EPJP was partially supported by the University of Iowa Department of Mathematics NSF VIGRE grant DMS-0602242.}

\author{Palle E. T. Jorgensen} 
\address{University of Iowa, Iowa City, IA 52246-1419 USA}
\ead{palle-jorgensen@uiowa.edu}

\author{Erin P. J. Pearse}
\address{University of Oklahoma, Norman, OK 73019-0315 USA}
\ead{ep@ou.edu}

%\date{\textbf{\today}}

\begin{abstract}
Let $(G,c)$ be an infinite network, and let $\mathcal{E}$ be the canonical energy form. Let $\Delta_2$ be the Laplace operator with dense domain in $\ell^2(G)$ and let $\Delta_{\mathcal{E}}$ be the Laplace operator with dense domain in the Hilbert space $\mathcal{H}_\mathcal{E}$ of finite energy functions on $G$. It is known that $\Delta_2$ is essentially self-adjoint, but that $\Delta_{\mathcal{E}}$ is \emph{not}. In this paper, we characterize the Friedrichs extension of $\Delta_{\mathcal{E}}$ in terms of $\Delta_2$ and show that the spectral measures of the two operators are mutually absolutely continuous with Radon-Nikodym derivative $\lambda$ (the spectral parameter),  in the complement of $\lambda=0$. We also give applications to the effective resistance on $(G,c)$.  
For transient networks, the Dirac measure at $\lambda = 0$ contributes to the spectral resolution of the Friedrichs extension of $\Delta_{\mathcal{E}}$ but not to that of the self-adjoint $\ell^2$ Laplacian.
%It turns out that it is the spectral theory of \Lap as an operator in \HE (as opposed to $\ell^2(G)$) which reflects important properties of the network, including certain precise notions of metric and boundary.
%  While there is a substantial literature (pure and applied) concerning graph-Laplacians on infinite networks, much less developed, in the use of Hilbert spaces over infinite graphs, is the distinction between operator in the $\ell^2$ space of the vertex-set $V$  vs the case when the Hilbert space is defined by a suitable energy form. By a network we mean a triple $(V, E, c)$ where $V$ is a (typically countable infinite) set of vertices in a graph, and with $E$ denoting the set of edges. The function $c$ in the triple is fixed and defined on $E$. It represents conductance in electrical network models, and assumed symmetric and positive on E. Starting with $(V, E, c)$, we introduce a graph-Laplacian $\Delta$ , and an energy Hilbert space \HE (both depending on $c$). While it is known that $\Delta$ is essentially selfadjoint on its natural dense domain in $\ell^2(V)$, its realization in \HE typically is not. We give a characterization of the Friedrichs extension of the \HE Laplacian, and we prove a formula for its computation. A number of spectral theoretical corollaries are derived from our main theorem.
\end{abstract}

%spectral permanence \sep singular equations, singular points, divergence form, time to infinity, 

\end{frontmatter}

%\linenumbers

\setcounter{tocdepth}{1} %{\footnotesize \tableofcontents}

\allowdisplaybreaks

%%!TEX root = Friedrichs.tex

\section{Introduction}

We study Laplace operators on infinite networks, and their self-adjoint extensions. Here, a network is just an connected undirected weighted graph $(G, c)$; see Definition~\ref{def:network}. The associated network Laplacian \Lap acts on functions $u:G \to \bR$; see Definition~\ref{def:graph-laplacian}. We restrict attention to the case when the network is \emph{transient}%
	\footnote{This equivalent to assuming the existence of \emph{monopoles}; see Definition~\ref{def:dipole} and Remark~\ref{rem:monotransience}.}, 
	and we are particularly interested in the case when \Lap is unbounded, in which case some care must be taken with the domains. We consider \Lap separately as an operator on \HE, the Hilbert space of finite energy functions on $G$ and as on operator on $\ell^2(G)$. 
	Although the two operators agree formally, their spectral theoretic properties are quite different.
    The space \HE is defined in terms of the quadratic form \energy, which gives the Dirichlet energy of a function $u$; see Definition~\ref{def:H_energy}. By $\ell^2(G)$, we mean the unweighted space of square-summable functions on $G$ under counting measure; see Definition~\ref{def:ell2}. 

Neither of the two Hilbert spaces is contained in the other, and the two Hilbert norms do not compare. It follows that, these two incarnations of the Laplacian may have quite different spectral theory. Common to the two is that \Lap is defined on its natural dense domain in each of the Hilbert spaces (these domains are given in Definition~\ref{def:domLapE} and Definition~\ref{def:domLap2}), and in each case it is a Hermitian and non-negative operator. 
However, it is known from \cite{Woj07, SRAMO, KellerLenz09, KellerLenz10} that \Lap is essentially self-adjoint on its natural domain in $\ell^2(G)$ but in \cite{SRAMO} it is shown that \Lap is \emph{not} essentially self-adjoint on its natural domain in \HE (see Definition~\ref{def:domLapE}). Nonetheless, we prove that the Friedrich extension of the latter has a spectral theory that can be compared with the former.

Theorem~\ref{thm:Friedrichs} is our first main result, and it characterizes the Friedrichs extension (see Definition~\ref{def:Friedrichs}) of the Laplacian on \HE in terms of the Laplacian on $\ell^2(G)$. Theorem~\ref{thm:spectral-resolutions} is our other main result, and it shows that the spectral measures of the Laplacian on $\ell^2(G)$ and the Friedrichs extension on \HE are mutually absolutely continuous with Radon-Nikodym derivative \gl (the spectral parameter). 
We use Theorem~\ref{thm:spectral-resolutions} to derive a number of spectral-theoretic conclusions. In particular, Corollary~\ref{thm:R(x,y)} gives a formula for the (effective) wired resistance metric on $(G,c)$ in terms of the spectral resolution of \Lap on $\ell^2(G)$, and Corollary~\ref{thm:R(x,y)-spectral-gap} shows that a spectral gap for \Lap on $\ell^2(G)$ implies a bound on the wired resistance. Resistance is a natural metric on networks and has been considered previously in many contexts; see \cite{Kig01,Kig03,Lyons,Soardi94,Woess09,DoSn84,Thomassen90}. Also see Definition~\ref{def:R(x)} and the ensuing discussion; some examples are explored in \S\ref{sec:examples}. We are also able to use spectral methods to recover some classical results for the integer lattices $(\bZd,\one)$ in \S\ref{sec:examples}. It turns out that it is the spectral theory of \Lap as an operator in \HE (as opposed to $\ell^2(G)$) which reflects important properties of the network, including certain precise notions of metric and boundary. The Friedrichs extension of the Laplacian arises naturally in this context as it corresponds to a limiting case of Dirichlet boundary conditions for the Laplacian, and hence to wired resistance metric; see \cite[Rem.~2.22]{ERM}. %it is the boundary of the network upon which the discrepancy between the free and wired resistance metrics depends.%; see \cite[Rem.~2.22]{ERM} for a discussion of resistance metrics in terms of boundary conditions.

%By contrast the discrete Laplacian \Lap to be studied here, the starting point is a prescribed network, typically an infinite graph with vertices $V$ and a conductance function $c$  defined on its edges $E$.   In our present study of the graph-Laplacian \Lap,  we stress an important distinction between (i) the operator theory for the $\ell^2$ space of the set $V$ of vertices, and (ii) the case when the Hilbert space \HE is defined by a suitable energy form (depending on c), the energy Hilbert space. 

To make this paper accessible to diverse audiences, we have included a number of definitions we shall need from the theory of (i) infinite networks, and (ii) the use of unbounded operators on Hilbert space in discrete contexts. Some useful background references for the first are \cite{Soardi94} and \cite{Woess09} and the multifarious references cited therein; see also \cite{Yamasaki79, Zem91, HK10, Kayano88, Kayano84, Kayano82, MuYaYo, vBL09, DuJo10}. For the second, see \cite{DuSc88} (especially Ch.~12) and \cite{vN32, Sto90, Jor80b, DJ06, BB09}.
%Connections between spectral theory and metric geometry has received attention in recent papers, see \cite{Chung, CaW92, Woess09} for example,  and in operator algebra theory \cite{Rie04, Rie99, Rie98}. 
For relevant background on reproducing kernels, see e.g., \cite{PaSc72, Aronszajn50, MuYaYo, Kal70}. %Kai65, AL08, Fuglede05
   In our first section below, we have recorded some lemmas from \cite{DGG, ERM, Multipliers, SRAMO, bdG, RBIN, RANR, Interpolation, LPS, OTERN} in the form in which they will be needed in the rest of the paper. Some of these results are folkloric or well known in the literature; in such cases, we refer to our own papers only for convenience.

\section{Basic terms and previous results}
\label{sec:Basic-terms-and-previous-results}

We now proceed to introduce the key notions used throughout this paper: resistance networks, the energy form \energy, the Laplace operator \Lap, and their elementary properties. 

\begin{defn}\label{def:network}
  A \emph{(resistance) network} $(\Graph,\cond)$ consists of a connected undirected graph \Graph and a symmetric \emph{conductance function} $\cond: \Graph \times \Graph \to [0,\iy)$. We write $x,y \in \Graph$ to indicate that $x$ and $y$ are vertices of the graph. The conductance function defines the adjacency relation as follows: $x$ and $y$ are neighbours (i.e., there is an edge connecting $x$ and $y$) iff $c_{xy}>0$, in which case the nonnegative number $c_{xy}=c_{yx}$ is the weight (conductance, or reciprocal resistance) associated to this edge. 
  
  We make the standing assumption that $(\Graph,\cond)$ is \emph{locally finite}. This means that every vertex has \emph{finite degree}, i.e., for any fixed $x \in \Graph$ there are only finitely many $y \in \Graph$ for which $c_{xy}>0$. We denote the net conductance at a vertex by 
  \linenopax
  \begin{align}\label{eqn:c(x)}
      \cond(x) := \sum_{y \nbr x} \cond_{xy}.     
  \end{align}
  %We require $\cond(x) < \iy$, but $\cond(x)$ need not be a bounded function on \Graph. The notation \cond may be used to indicate the multiplication operator $(\cond v)(x) := \cond(x) v(x)$.
  %i.e., the diagonal matrix with entries $\cond(x)$ with respect to the (vector space) basis $\{\gd_x\}$.

In this paper, \emph{connected} means simply that for any $x,y \in \Graph $, there is a finite sequence $\{x_i\}_{i=0}^n$ with $x=x_0$, $y=x_n$, and $\cond_{x_{i-1} x_i} > 0$, $i=1,\dots,n$.  
%We may assume there is at most one edge from $x$ to $y$, as two conductors $\cond^1_{xy}$ and $\cond^2_{xy}$ connected in parallel can be replaced by a single conductor with conductance $\cond_{xy} = \cond^1_{xy} + \cond^2_{xy}$. Also, we assume $\cond_{xx}=0$ so that no vertex has a loop. 

%Since the edge data of $(\Graph,\cond)$ is carried by the conductance function, we will henceforth simplify notation and write $x \in \Graph$ to indicate that $x$ is a vertex. 
For any network, one can fix a reference vertex, which we shall denote by $o$ (for ``origin''). It will always be apparent that our calculations depend in no way on the choice of $o$.
\end{defn}

\begin{defn}\label{def:graph-laplacian}
  The \emph{Laplacian} on \Graph is the linear difference operator which acts on a function $u:\Graph \to \bR$ by
  \linenopax
  \begin{equation}\label{eqn:Lap}
    (\Lap u)(x) :
    = \sum_{y \nbr x} \cond_{xy}(u(x)-u(y)).
  \end{equation}
  A function $u:\Graph \to \bR$ is \emph{harmonic} iff $\Lap u(x)=0$ for each $x \in \Graph$. 
\end{defn}

The domain of \Lap, considered as an operator on \HE or $\ell^2(G)$, is given in Definition~\ref{def:domLapE} and Definition~\ref{def:domLap2}.

\begin{comment}[Conventions for \Lap]
  \label{rem:Lap-conventions}
  We have adopted the physics convention%
  \footnote{Our motivation for this convention is that we prefer to have the spectrum of \Lap be nonnegative. This will hold for \Lap as an operator on \HE and on $\ell^2$, according the definitions (domains) given below. See Remark~\ref{rem:ell2(c)}.}  
  and thus our Laplacian is the negative of the one commonly found in the PDE literature. The network Laplacian \eqref{eqn:Lap} should not be confused with the stochastically renormalized Laplace operator $\cond^{-1} \Lap$ which appears in the probability literature, or with the spectrally renormalized Laplace operator $\cond^{-1/2} \Lap \cond^{-1/2}$ which appears in the literature on spectral graph theory (e.g., \cite{Chung}). Here, $c$ denotes the multiplication operator associated to the function $c(x)$. See Remark~\ref{rem:ell2(c)}.
\end{comment}

\subsection{The energy space \HE} 
\label{sec:The-energy-space}

\begin{defn}\label{def:graph-energy}
  The \emph{energy form} is the (closed, bilinear) Dirichlet form
  \linenopax
  \begin{align}\label{eqn:def:energy-form}
    \energy(u,v)
    := \frac12 \sum_{x,y \in \Graph} \cond_{xy}(u(x)-u(y))(v(x)-v(y)),
  \end{align}
  which is defined whenever the functions $u$ and $v$ lie in the domain
  \linenopax
  \begin{equation}\label{eqn:def:energy-domain}
    \dom \energy = \{u:\Graph \to \bR \suth \energy(u,u)<\iy\}.
  \end{equation}
  Hereafter, we write the energy of $u$ as $\energy(u) := \energy(u,u)$. Note that $\energy(u)$ is a sum of nonnegative terms and hence converges iff it converges absolutely. %Since $\cond_{xy}=\cond_{yx}$ and $\cond_{xy}=0$ for nonadjacent vertices, the initial factor of $\frac12$ in \eqref{eqn:def:energy-form} corresponds to the fact that there is exactly one term in the sum for each edge in the network. 
\end{defn}

\begin{comment}\label{rem:sum-convergence}
  To remove any ambiguity about the precise sense in which \eqref{eqn:def:energy-form} converges, note that $\energy(u)$ is a sum of nonnegative terms and hence converges iff it converges absolutely. Since the Schwarz inequality gives $\energy(u,v)^2 \leq \energy(u)\energy(v)$, it is clear that the sum in \eqref{eqn:def:energy-form} is well-defined whenever $u,v \in \dom \energy$.
\end{comment}

\begin{comment}[\bR vs. \bC]
  \label{rem:RvsC}
  %For most of the analysis we carry out, i
  It will suffice to consider only \bR-valued functions for now. 
  %In particular, Lemma~\ref{thm:energy-kernel-is-real-and-bounded} (and the fact that the Laplacian commutes with conjugation) implies that m
  Most results extend immediately to the \bC-valued case with only trivial modification, for example, \eqref{eqn:def:energy-form} becomes $$\energy(u,v) = \frac12 \sum_{x,y \in \Graph} \cond_{xy}(\cj{u(x)}-\cj{u(y)})(v(x)-v(y)).$$ %For further details, please see \cite{bdG}.
\end{comment}

The energy form \energy is sesquilinear and conjugate symmetric on $\dom \energy$ and would be an inner product if it were positive definite. Let \one denote the constant function with value 1 and observe that $\ker \energy = \bR \one$. One can show that $\dom \energy / \bR \one$ is complete and that \energy is closed;  see \cite{DGG,OTERN}, \cite{Kat95}, or \cite{FOT94}.

\begin{defn}\label{def:H_energy}\label{def:The-energy-Hilbert-space}
  The \emph{energy (Hilbert) space} is $\HE := \dom \energy / \bR \one$. The inner product and corresponding norm are denoted by
  \linenopax
  \begin{equation}\label{eqn:energy-inner-product}
    \la u, v \ra_\energy := \energy(u,v)
    \q\text{and}\q
    \|u\|_\energy := \energy(u,u)^{1/2}.
  \end{equation}
\end{defn}

It is shown in \cite[Lem.~2.5]{DGG} that the evaluation functionals $L_x u = u(x) - u(o)$ are continuous, and hence correspond to elements of \HE by Riesz duality (see also \cite[Cor.~2.6]{DGG}).

\begin{defn}\label{def:vx}\label{def:energy-kernel}
  Let $v_x$ be defined to be the unique element of \HE for which
  \linenopax
  \begin{equation}\label{eqn:v_x}
    \la v_x, u\ra_\energy = u(x)-u(o),
    \qq \text{for every } u \in \HE.
  \end{equation}
  Note that $v_o$ corresponds to a constant function, since $\la v_o, u\ra_\energy = 0$ for every $u \in \HE$. Therefore, $v_o$ may be safely omitted in some calculations. 
\end{defn}

  Equation \eqref{eqn:v_x} means that the collection $\{v_x\}_{x \in \Graph}$ forms a reproducing kernel for \HE and thus has dense span in \HE. We call $\{v_x\}_{x \in \Graph}$ the \emph{energy kernel}.

\begin{remark}[Differences and representatives]\label{rem:differences}
  Equation \eqref{eqn:v_x} is independent of the choice of representative of $u$ because the right-hand side is a difference: if $u$ and $u'$ are both representatives of the same element of \HE, then $u' = u+k$ for some $k \in \bR$ and
  $u'(x) - u'(o) = (u(x)+k)-(u(o)+k) = u(x)-u(o).$
  By the same token, the formula for \Lap given in \eqref{eqn:Lap} describes unambiguously the action of \Lap on equivalence classes $u \in \HE$. Indeed, formula \eqref{eqn:Lap} defines a function $\Lap u:\Graph \to \bR$ but we may also interpret $\Lap u$ as the class containing this representative.
\end{remark}

\begin{defn}\label{def:d_x}
  Let $\gd_x \in \ell^2(G)$ denote the Dirac mass at $x$, i.e., the characteristic function of the singleton $\{x\}$ and let $\gd_x \in \HE$ denote the element of \HE which has $\gd_x \in \ell^2(G)$ as a representative. The context will make it clear which meaning is intended. Observe that $\energy(\gd_x) = \cond(x) < \iy$ is immediate from \eqref{eqn:def:energy-form}, and hence one always has $\gd_x \in \HE$ (recall that $c(x)$ is the total conductance at $x$; see \eqref{eqn:c(x)}).
\end{defn}

\begin{defn}\label{def:Fin}
  For $v \in \HE$, one says that $v$ has \emph{finite support} iff there is a finite set $F \ci G$ such that $v(x) = k \in \bC$ for all $x \notin F$. Equivalently, the set of functions of finite support in \HE is 
  \linenopax
  \begin{equation}\label{eqn:span(dx)}
    \spn\{\gd_x\} = \{u \in \dom \energy \suth u(x)=k \text{ for all } x \notin F\},
    %some $k$, for all but finitely many } x \in G\},
  \end{equation}
  for some finite $F \ci G$. 
  Define \Fin to be the \energy-closure of $\spn\{\gd_x\}$. 
\end{defn}

\begin{defn}\label{def:Harm}
  The set of harmonic functions of finite energy is denoted
  \linenopax
  \begin{equation}\label{eqn:Harm}
    \Harm := \{v \in \HE \suth \Lap v(x) = 0, \text{ for all } x \in G\}.
  \end{equation}
\end{defn}

The following result is well known; see \cite[\S{VI}]{Soardi94}, \cite[\S9.3]{Lyons}, \cite[Thm.~2.15]{DGG}, or the original \cite[Thm.~4.1]{Yamasaki79}.

\begin{theorem}[Royden Decomposition]\label{thm:HE=Fin+Harm}
  $\HE = \Fin \oplus \Harm$.
\end{theorem}
 
%Theorem~\ref{thm:HE=Fin+Harm} first appeared in \cite[Thm.~4.1]{Yamasaki79}, where it was called the ``Royden Decomposition'' by analogy with Royden's result for Riemann surfaces. However, the result is incorrect as stated there and the present corrected form may be found in \cite[\S{VI}]{Soardi94} or \cite[\S9.3]{Lyons}. This result also follows immediately from Lemma~\ref{thm:<delta_x,v>=Lapv(x)}; see \cite[Thm.~2.15]{DGG}. %In view of Theorem~\ref{thm:HE=Fin+Harm}, denote the (orthogonal) projection to \Fin by \Pfin and the projection to \Harm by \Phar.

\begin{defn}\label{def:dipole}
  A \emph{monopole} is any $w \in \HE$ satisfying the pointwise identity $\Lap w = \gd_x$ (in either sense of Remark~\ref{rem:differences}) for some vertex $x \in \Graph$. 
  A \emph{dipole} is any $v \in \HE$ satisfying the pointwise identity $\Lap v = \gd_x - \gd_y$ for some $x,y \in \Graph$.%
  %\footnote{Here, $\gd_x$ can be interpreted as either a function on $G$ or an element of \HE; both have the same meaning. See Remark~\ref{rem:differences} and Definition~\ref{def:d_x}.} 
\end{defn}

\begin{remark}\label{rem:monotransience}
  It is easy to see from the definitions (or \cite[Lemma~2.13]{DGG}) that energy kernel elements are dipoles, i.e., that $\Lap v_x = \gd_x - \gd_o$, and that one can therefore always find a dipole for any given pair of vertices $x,y \in G$, namely, $v_x-v_y$. On the other hand, monopoles exist if and only if the network is transient (see \cite[Thm.~2.12]{Woess00} or \cite[Rem.~3.5]{DGG}). %, which is the case in which we are primarily interested. This holds, for example, whenever $\Harm \neq 0$ (see \cite[Cor.~4.3]{DGG}). Note that if $w_x$ is a monopole at $x \in \Graph$, then $w-v_x+v_y$ is a monopole at $y \in \Graph$, and so there is a monopole at $x$ iff there is a monopole at every vertex. This famous fact is better known in probabilistic language: the random walk started at $x$ is transient iff it is transient when started from any vertex. See, e.g. \cite[\S1--2]{Woess00}.
\end{remark}  

\begin{comment}\label{thm:minimal-monopole}
  On any transient network, there is a unique energy-minimizing monopole at each vertex.  
  \begin{proof}
    Any two monopoles at $x$ must differ by an element of \Harm, so write $\tilde w = \Pfin w$ for any monopole $w$ at $x$. Then $\tilde w$ is energy-minimizing because $\energy(w) = \energy(\tilde w) + \energy(\Phar w)$ by Pythagoras' theorem.
  \end{proof}
\end{comment}

%Alternatively, Lemma~\ref{thm:minimal-monopole} holds because the space of monopoles at $x$ is closed and convex, so the quadratic form \energy has a unique minimum here.
  
\begin{remark}\label{rem:normalization}
  Denote the unique energy-minimizing monopole at $o$ by $w_o$; the existence of such an object is explained in \cite[\S3.1]{DGG}.
  We will be interested in the family of monopoles defined by
  \linenopax
  \begin{align}\label{eqn:w_x}
    \monov := w_o + v_x, \qq x \neq o.
  \end{align}
  In \S\ref{sec:Phi} (see Lemmas~\ref{thm:Lap-ripple}--\ref{thm:Lap-under}) we use the representatives specified by
  %\footnote{With the appropriate definition of $R(x,\iy)$ as the effective resistance between $x$ and ``the set at \iy'', it turns out that $w_x(x) = \energy(w_x) = R(x,\iy)$; see Definition~\ref{def:R(x)}. More precisely, put $y=\iy_k$ in \cite[(2.20)]{ERM};}
  \linenopax
  \begin{align}\label{eqn:w_x(o)}
    \monov(x) = \energy(\monov), % = R(x,\iy), 
    \qq\text{and}\qq
    v_x(o)=0.
  \end{align}
  When $\Harm=0$, $\energy(\monov)$ is the \emph{capacity} of $x$; see, e.g., \cite[\S4.D]{Woess09}.
  %See \cite{ERM} for more about $R(x,\iy)$. 
  %In this case, one representative of $w_o$ is given by $w_o(x) = g(x,o)/c(o)$ where $g(x,y)$ is the Green kernel on $(G,c)$ defined in terms of the transition probabilities $p(x,y) := c_{xy}/c(x)$; see \cite[Rem.~3.3]{DGG} for further discussion. It follows from \eqref{eqn:w_x(o)} and the symmetry $g(x,o)/c(o) = g(o,x)/c(x)$ that 
  %\linenopax
  %\begin{align}\label{eqn:monopole-symmetry}
  %  w_x(y) = w_y(x), \qq \text{for all } x,y \in G.
  %\end{align}
\end{remark}

\begin{lemma}[{\cite[Lem.~2.11]{DGG}}]
  \label{thm:<delta_x,v>=Lapv(x)}
  For $x \in \Graph$ and $u \in \HE$,  $\la \gd_x, u \ra_\energy = \Lap u(x)$.
  \begin{proof}
    Compute $\la \gd_x, u \ra_\energy = \energy(\gd_x, u)$ directly from formula \eqref{eqn:def:energy-form}.
  \end{proof}
\end{lemma}

\begin{lemma}\label{thm:Lap-mono-Kron}
  For any $x,y \in G$, 
  \linenopax
  \begin{align}\label{eqn:Lap-mono-Kron}
    \Lap \monov(y) = \Lap \monoy(x) 
    = \la \monov,\Lap \monoy\ra_\energy 
    = \la \Lap \monov, \monoy\ra_\energy = \gd_{xy},  
  \end{align}
  where $\gd_{xy}$ is the Kronecker delta.
  \begin{proof}
    First, note that $\Lap \monov(y) = \gd_{xy} = \Lap \monov(y)$ as functions, immediately from the definition of monopole. Then Lemma~\ref{thm:<delta_x,v>=Lapv(x)} gives $\la \monov,\Lap \monoy\ra_\energy
      = \la \monov,\gd_y\ra_\energy
      = \Lap \monov(y)$    since $\Lap \monoy = \gd_y$ and $\la u,\gd_y\ra_\energy = \Lap u(y)$, and similarly for the other identity. 
  \end{proof}
\end{lemma}

\begin{defn}\label{def:domLapE}
  On \HE, start with \Lap defined on $\spn\{\monov\}_{x \in G}$ pointwise by \eqref{eqn:Lap}, %the subspace of (finite) linear combinations of monopoles, 
  and then obtain the closed operator \LapE by taking the graph closure; the following lemma shows that this is justified.
  %to be the symmetric operator which is given pointwise by \eqref{eqn:Lap} and which has dense domain
  %\linenopax
  %\begin{align}\label{eqn:domLapE}
  %  \dom \LapE := \{u \in \HE \suth \exists (u_n)_{n=1}^\iy \ci \spn\{w_x\}_{x \in G} \text{ with } \lim_{n \to \iy}\||u-u_n\||=0\}, 
  %\end{align}
  %where $\||\cdot\||$ denotes the graph norm. 
  %
  %In other words, $u \in \dom \LapE$ if and only if there is a sequence $(u_n)_{n=1}^\iy \ci \spn\{w_x\}_{x \in G}$ and a vector $v \in \HE$ for which 
  %\linenopax
  %\begin{align*}%\label{eqn:}
  %  \lim_{n \to \iy} \|u-u_n\|_\energy = 0
  %  \qq\text{and}\qq
  %  \lim_{n \to \iy} \|v-\Lap u_n\|_\energy = 0,
  %\end{align*}
  %in which case one writes $\LapE u = v$.
\end{defn}

%A closely related form of the following fact first appeared in \cite[Lem.~3.15]{SRAMO}.
 
\begin{lemma}\label{thm:semibounded}
   \LapE is well-defined and non-negative (hence also closed and Hermitian).%. A fortiori, \LapE is . 
  \begin{proof}
    Let $\gx \in \dom \LapE$ with $\spt \gx$ contained in some finite set $F \ci G$. By \eqref{eqn:Lap-mono-Kron}, 
    %Computing directly from the left side of \eqref{eqn:semibounded}, one has 
    \linenopax
    \begin{align}\label{eqn:semibounded}
      \la u, \Lap u \ra_\energy
      &= \sum_{x,y \in F} {\gx_x} \gx_y \la \monov,\Lap \monoy\ra_\energy 
       = \sum_{x,y \in F} {\gx_x} \gx_y \gd_{xy}
       = \sum_{x \in F} |\gx_x|^2 \geq 0.  
    \end{align} 
    The closure of any semibounded operator is semibounded. This implies \LapE is Hermitian and hence contained in its adjoint. Since every adjoint operator is closed, \LapE is closable.
  \end{proof}
\end{lemma}

%\begin{remark}\label{rem:semibounded}
%  Typically, an operator $T$ is called \emph{semibounded} iff there is a $b\in \bR$ such that either $\la u, Tu\ra \geq b$ \textbf{or} $\la u, Tu\ra \leq b$ holds for all $u \in \dom T$. $T$ is a \emph{non-negative} operator in the special case when $\la u, Tu\ra \geq 0$ for all $u \in \dom T$. Although the specific operators 
%\end{remark}

\subsection{The Hilbert space $\ell^2(G)$} 
\label{sec:ell2(G)}

As there are many uses of the notation $\ell^2(G)$, we provide the following elementary definitions to clarify our conventions.
%  Note that when considering $\ell^2(G)$, we use 

\begin{defn}\label{def:ell2}
For functions $u,v:G \to \bR$, define the inner product
  \linenopax
  \begin{align}\label{eqn:ell2-inner-product}
    \la u, v\ra_2 := \sum_{x \in G} u(x) v(x).
  \end{align}
  %i.e., the measure on $\ell^2(G)$ is simply counting measure. %Note that \eqref{eqn:ell2-inner-product} is written for \bR-valued functions $u$ and $v$; see Remark~\ref{rem:RvsC}.
\end{defn}
  
\begin{defn}\label{def:domLap2}
  On $\ell^2(G)$, we begin with \Lap defined pointwise by \eqref{eqn:Lap} on $\spn\{\gd_x\}_{x \in G}$, the subspace of (finite) linear combinations of point masses, and then obtain the closed operator $\Lap_2$ by taking the graph closure (see Remark~\ref{rem:domLap2}).
\end{defn}

\begin{remark}\label{rem:domLap2}
  \cite[Lem.~2.7 and Thm.~2.8]{SRAMO} states that $\Lap_2$ is semibounded and essentially self-adjoint. It follows that $\Lap_2$ is closable by the same arguments as in the end of the proof of Lemma~\ref{thm:semibounded}.
  See also \cite{Woj07, KellerLenz09, KellerLenz10}.
\end{remark}

%%!TEX root = Friedrichs.tex

\section{Some properties of the Laplacian and the monopoles}

\begin{lemma}\label{thm:Lap-ripple}
  For any $x \in G$, 
  \linenopax
  \begin{align}\label{eqn:Lap-ripple}
    \gd_x = c(x) \monov - \sum_{y \nbr x} c_{xy} \monoy.  
  \end{align}
  \begin{proof}
    For any $z \in G$, formulas \eqref{eqn:w_x(o)}, \eqref{eqn:Lap-mono-Kron} and \eqref{eqn:Lap} give
    \linenopax
    \begin{align*}%\label{eqn:}
      \gd_x(z) = \gd_z(x)
      &= \Lap \monov[z](x)
      = c(x)\monov[z](x) - \sum_{y \nbr x} c_{xy} \monov[z](y)
      = c(x)\monov(z) - \sum_{y \nbr x} c_{xy} \monoy(z).
      \qedhere
    \end{align*} 
  \end{proof}
\end{lemma}

\begin{lemma}\label{thm:Lap-under}
  For any $x,y \in G$, if $\gd_{xy}$ is the Kronecker delta and $\Lap_x$ denotes the Laplacian taken with respect to the $x$ variable, then
  \linenopax
  \begin{align}\label{eqn:Lap-under}
    \Lap_x \la \monov, \monoy\ra_\energy = \la \Lap \monov, \monoy\ra_\energy = \gd_{xy}.  
  \end{align}
  \begin{proof}
    Using \eqref{eqn:Lap}, we have
    \linenopax
    \begin{align*}%\label{eqn:}
      \Lap_x \la \monov, \monoy\ra_\energy
      = c(x) \la \monov, \monoy\ra_\energy - \sum_{z \nbr x} c_{xz} \la \monov[z], \monoy\ra_\energy
      = \left\la c(x) \monov - \sum_{z \nbr x} c_{xz} \monov[z], \monoy\right\ra\SubEn
    \end{align*} 
    whence the result follows by applying \eqref{eqn:Lap-ripple} and then \eqref{eqn:Lap-mono-Kron}. 
  \end{proof}
\end{lemma}

\subsection{The transformation $\gF:\gd_x \mapsto \monov$ and the matrix $M=[\la \monov, \monov[y]\ra_\energy]_{x,y \in G}$}
\label{sec:Phi}

\begin{defn}\label{def:Phi}
  %For $\gx \in \sB$, 
  Define $\gF : \ell^2(G) \to \HE$ on $\dom\gF = \spn\{\gd_x\}_{x \in G}$ by $\gF\gd_x = \monov$. 
\end{defn}

  Note that $\ran \gF$ is dense in \HE because it contains $\spn\{v_x\}_{x \in G}$; see \eqref{eqn:w_x}. % and Definition~\ref{def:energy-kernel}. 

\begin{remark}\label{rem:closable}
  The operator \gF may not be closable.
  %; equivalently, the adjoint may not be densely defined. As an example, if $\monov \in \dom \gF^\ad$,%
  %\footnote{This happens, for instance, on the binary tree example discussed in \S\ref{sec:tree2}.} then
  %\linenopax
  %\begin{align*}%\label{eqn:}
  %  (\gF^\ad \monoy)(x)
  %  = \la \gd_x, \gF^\ad \monoy\ra_{\ell^2}
  %  = \la \gF \gd_x, \monoy\ra_{\energy} 
  %  = \la \monov, \monoy\ra_{\energy}.
  %\end{align*}
  %Since this is true for every $x,y \in G$, it must be that $\gF^\ad \monoy$ is the representative of \monoy specified in \eqref{eqn:w_x(o)}. Since $\dom\gF^\ad$ 
  %$\{u \in \HE \suth \gF^\ad u \in \ell^2(G)\}$ 
  %is a subspace of \Fin and hence cannot be dense in \HE whenever $\Harm \neq 0$; see \cite{DGG} for examples. 
  This necessitates some care in the formulation of the Friedrichs extension in Definition~\ref{def:Phi-ext}. 
\end{remark}

\begin{defn}\label{def:M}
  %For a finite set $F \ci G \less \{o\}$, we have $\gF\gx = \sum_{x \in F} \gx(x) \monov$, for all $\gx \in \spn\{\gd_x\}$. 
  Let $M$ be the (infinite) matrix with entries 
  $M_{xy} = \la \monov, \monoy \ra_{\energy}$. 
  %and let $M_F := M\evald{F \times F}$ be the submatrix of $M$ defined by deleting all rows and columns corresponding to points $x \notin F$, i.e., $M_F$ is an $|F| \times |F|$ matrix with entries
  %\linenopax
  %\begin{equation}\label{eqn:def:MF}
  %  (M_F)_{xy} = \la \monov,\monoy \ra_\energy, 
  %  \q \forall x,y \in F.
  %\end{equation}
  %In general, one may have $\monov \in \spn\{\monov\}_{x \in F}$ with support extending outside of $F$; examples are given in \cite{OTERN}. 
\end{defn}

\begin{lemma}\label{thm:M-intertwines}
  For all $\gx \in \dom \gF$, one has $\la\gx, M\gx \ra_{\ell^2} = \| \gF(\gx)\|_\energy^2$.
  \begin{proof}  
    The computation is immediate: 
    \linenopax
    \begin{align*}%\label{eqn:}
      \| \gF(\gx)\|_\energy^2
      & = \left\la \sum_{x \in F} \gx(x) \monov, \sum_{y \in F} \gx(y) \monoy \right\ra \SubEn
      = \sum_{x \in F} \sum_{y \in F} {\gx(x)} \gx(y) M_{x,y}
      = \la \gx, M \gx\ra_{\ell^2}.
      \qedhere
    \end{align*}
  \end{proof}
\end{lemma}

Lemma~\ref{thm:M-intertwines} shows that the matrix $M$ plays the role of the formal expression $\gF^\ad\gF$.
We will need a couple of lemmas relating $\Lap_2$ and \LapE to \gF. Lemma~\ref{thm:Phi-intertwines} relates the inner product of \HE to the inner product on $\ell^2(G)$, and shows that they differ ``by a Laplacian''; see also \cite[Lem.~5.30]{SRAMO}. %Lemma~\ref{thm:intertwined-powers} describes how \gF intertwines iterates of \Lap. In these lemmas, all the functions have finite support and so we suppress reference to domains for simplicity.

\begin{lemma}\label{thm:Phi-intertwines}
  For all $\gx,\gh \in \dom\gF$, one has $\la \gF\gx, \LapE \gF\gh \ra_\energy = \la\gx,\gh\ra_2.$
  \begin{proof}  
    Letting $\gx = \sum_{x \in F} \gx_x \gd_x$ and $\gh = \sum_{x \in F} \gh_x \gd_x$ (for some finite $F \ci G$) and arguing as in the proof of Lemma~\ref{thm:semibounded}, we have
    \linenopax
    \begin{align*}%\label{eqn:Phi-intertwines-der1}
      \la \gF\gx, \LapE \gF\gh \ra_\energy
      &= \sum_{x,y \in F} {\gx_x} \gh_y \la \monov,\LapE \monoy\ra_\energy 
       = \sum_{x,y \in F} {\gx_x} \gh_y \gd_{xy}
       = \la\gx,\gh\ra_2.
       \qedhere
    \end{align*} 
  \end{proof}
\end{lemma}

\begin{cor}\label{thm:intertwined-powers}
  For $\gx \in \dom\gF$, one has $\LapE\gF \gx = \gF \Lap_2\gx = \sum_{x \in F} \gx_x \gd_x$.
  %\linenopax
  %\begin{align}\label{eqn:intertwined-powers}
  %  \Lap^n \gF \gx = \gF \Lap^n \gx, % = \sum_{x \in F} \gx_x \gd_x,
  %  \qq\text{ for } n=0,1,2,\dots.
  %\end{align}
  \begin{comment}
    %\marginpar{Can we remove "banded"? This appears to be the only place where it is used.}
    First, note that $\Lap\gx$ has finite support because \Lap is banded (see Remark~\ref{rem:banded}) and \gx has finite support. Consequently, $\Lap\gx \in \dom\gF$. % because
    %\linenopax
    %\begin{align*}%\label{eqn:}
    %  \sum_{x \in F} \Lap\gx(x) 
    %  = \sum_{x \in F} \sum_{y \in F} c_{xy} \gx(y)  
    %  = \sum_{y \in F} \gx(y) \sum_{x \in F} c_{xy} 
    %  = 0.
    %\end{align*}
    On the one hand, % has $\gx = \sum_{x \in F} \gx_x \gd_x$
    \linenopax
    \begin{align}\label{eqn:intertwined-powers-der1}
      \Lap \gF \gx
      &= \Lap \sum_{x \in F} \gx_x \monov
       = \sum_{x \in F} \gx_x \Lap \monov
       = \sum_{x \in F} \gx_x \gd_x      
    \end{align}
    %since $\sum_{x \in F} \gx_x = 0$ for $\gx \in \dom\gF$. 
    by Lemma~\ref{thm:Lap-mono-Kron}. On the other hand, 
    \linenopax
    \begin{align*}%\label{eqn:Lap-entries}
      \gF \Lap \gd_x
      &= \gF \left(c(x) \gd_x - \sum_{y \in F} c_{xy} \gd_y\right)
       = c(x) \monov - \sum_{y \in F} c_{xy} \monoy
       %=\Lap \monov
       =\gd_x.   
    \end{align*}
    by Lemma~\ref{thm:Lap-ripple}, and $\gF \Lap \gx = \sum_{x \in F} \gx_x \gd_x$ follows by taking linear combinations. This establishes that $\Lap \gF \gx = \gF \Lap \gx$, whence \eqref{eqn:intertwined-powers} follows by iteration.
  \end{comment}
\end{cor}

Note that Corollary~\ref{thm:intertwined-powers} is an identity in \HE, and for this reason we have written $\sum_{x \in F} \gx_x \gd_x$ and not \gx, so as to account for the possible projection to \Fin.

%%!TEX root = Friedrichs.tex

\section{Characterization of the Friedrichs extension}
%\section{Friedrich's extension}

It is known from \cite[Prop.~4.9]{SRAMO} that \LapE may fail to be essentially self-adjoint. Therefore, we construct a canonical self-adjoint extension of \LapE, following the methods of Friedrichs (and von Neumann); see \cite[\S{}XII.5]{DuSc88} for background.

\begin{remark}%[$\Lap_2$ is essentially self-adjoint, but \LapE may not be]
  \label{rem:essentially-self-adjoint}
  Recall from Remark~\ref{rem:domLap2} that $\Lap_2$ is essentially self-adjoint and hence has a unique and well-defined spectral representation. \emph{We henceforth assume  without loss of generality that $\Lap_2$ is self-adjoint.}
  %Our task is to describe the Friedrichs extension of \LapE.   
%We also discuss $\Lap_2^{-1/2}$, but we define its domain in terms of the domain of a quadratic form (see \cite{Kat95}), rather than directly in terms of the functional calculus.
\end{remark}

In this section, we relate the domains given in Definition~\ref{def:domLapE} and Definition~\ref{def:domLap2}. This will entail comparing a quadratic form $q_M$ from $\ell^2(G)$ with a quadratic form $q_\Lap$ from \HE. We will have occasion to use the following result, which is a special case of {\cite[Ch.~VI, Thm.~2.1 and Thm.~2.6]{Kat95}}.

\begin{theorem}[Kato's Theorem]\label{thm:Kato-rep}
  Let $q$ be a densely defined, closed, symmetric sesquilinear form in a Hilbert space \sH which satisfies $\inf\{q(u) \suth u \in \dom q\} = 0$, for $q(u)=q(u,u)$.
  Then there is a unique self-adjoint operator $T$ with $\inf\{\la u, Tu\ra \suth u \in \dom T\} = \inf \spec T = 0$
  %which is (with lower bound $b$) 
  satisfying
  \begin{enumerate}[(i)]
    \item $\dom T \ci \dom T^{1/2} = \dom q$.
    \item  $q(u,v) = \la Tu, v\ra$ for any $u \in \dom T$ and $v \in \dom q$. 
    \item If $u \in \dom q$, $w \in \sH$, and $q(u,v) = \la w, v\ra$ for every $v \in \dom T$, \\ \hstr[6] then $u \in \dom T$ and $Tu=w$.
  \end{enumerate}
\end{theorem}

The next step is to extend the mapping $\Phi:\gd_x \to \monov$ from Definition~\ref{def:Phi} to functions which do not have finite support in Definition~\ref{def:Phi-ext}; this requires some further development of $M$ from  Definition~\ref{def:M}. %The extension of \gF appears in Lemma~\ref{thm:D} and Remark~\ref{rem:3-properties-of-completion}.

\begin{defn}\label{def:qM}
  The real Hermitian (symmetric) matrix $M$ defines a quadratic form with dense domain in $\ell^2(G)$. Define $q_M$ to be the closure of this form; note that this is justified by Kato's theorem because Lemma~\ref{thm:M-intertwines} shows that $M$ is non-negative.
\end{defn}

\begin{defn}\label{def:sM}
  %Consider the operation of multiplication by $M$, on $\spn\{\gd_x\}$. 
  Let \sM denote the self-adjoint operator corresponding to the quadratic form $q_M$ by Kato's theorem. 
\end{defn}

%By abuse of notation, we let $M$ denote the multiplication operator induced by the matrix $M$ in the following lemma.

\begin{remark}\label{rem:HE-intertwine-vs-L2-intertwine}
  Lemma~\ref{thm:Green} is a renormalized version (or a symmetrized version; see Remark~\ref{rem:normalization}) of the standard identity that the Laplacian and Green operator are inverses. In this context, the proof comes by comparing quadratic forms associated to $M$ and to \Lap. In a different context, the question of comparing two quadratic forms, and closability, comes up in the study of Gaussian stochastic processes, see \cite[\S5]{AJL11} and \cite{AlJo11}.
\end{remark}

\begin{lemma}\label{thm:Green}
  $\Lap_2$ and \sM are inverses of each other:
  \begin{enumerate}[(i)]
  \item $\sM \Lap_2 \gh = \gh$, for any $\gh \in \dom \Lap_2$,  
   \q\text{and}\q (ii) $\Lap_2 \sM \gx = \gx$, for any $\gx \in \dom \sM$. 
  \end{enumerate}
  %$M = \Lap_2^{-1}$.
  \begin{proof}
    Since $\Lap_2$ is a self-adjoint operator in $\ell^2(G)$ and $\spn\{\gd_x\}$ is contained in $\dom \Lap_2^n$ for any $n\geq1$, the matrix of $\Lap_2$ relative to the onb $\spn\{\gd_x\}$ is 
    \linenopax
    \begin{align*}%\label{eqn:}
      \tilde \Lap_{x,y}
      := \la \gd_x, \Lap_2 \gd_y\ra_2
       = \begin{cases}
         c(x), & y=x,\\ -c_{xy}, & y \nbr x, \\ 0, &\text{else}.
       \end{cases}
    \end{align*}
    The following matrix multiplication uses Lemma~\ref{thm:Lap-mono-Kron} to show that $\tilde \Lap M = \id$:
    \linenopax
    \begin{align}\label{eqn:LapM=I}
      (\tilde \Lap M)_{x,y}
      &= \sum_{z \in X} \tilde\Lap_{x,z} M_{z,y}
       = \sum_{z \in X} \tilde\Lap_{x,z} \la \monov[z], \monov[y]\ra_\energy
       = \Lap_x \la \monov[x], \monov[y]\ra_\energy
       = \gd_{xy}.
    \end{align}
    Note that the summation over $z$ is finite because $\tilde \Lap$ is banded; see Remark~\ref{rem:banded}. The computation for $M \tilde \Lap = \id$ is identical by the symmetry of the matrices. Therefore, $\Lap_2$ and \sM are inverses on a formal level.
    %It follows that that $\la \gd_x, \tilde \Lap M \gd_y\ra_2 = \gd_{xy}$.
    %The foregoing computation shows , and hence $M$ maps $\spn\{\gd_x\}$ into $\dom \Lap_2$, and that $\tilde \Lap M = \id$. To see this, 
    
    For (i), let $\gh \in \dom \Lap_2$. 
    % and define $\gx := \Lap_2\gh$, i.e., suppose that $(\gh_n)_{n=1}^\iy \ci \spn\{\gd_x\}$ and that there exists $\gx \in \ell^2(G)$ satisfying 
    %\linenopax
    %\begin{align*}%\label{eqn:}
    %  \lim_{k \to \iy} \|\gh - \gh_k\|_2 = 0
    %  \qq\text{and}\qq \lim_{k \to \iy} \|\gx - \Lap_2 \gh_k\|_2 = 0.
    %\end{align*}
    Since $\Lap_2$ is banded, the following double sum is finite:
    %$\Lap_2 \gh_k \in \spn\{\gd_x\} \ci \dom \sM$ for every $k$, and
    \linenopax
    \begin{align*}%\label{eqn:}
      q_M(\Lap_2 (\gh_n - \gh_m))
      &= \sum_{x,y} (\tilde \Lap (\gh_n - \gh_m))(x) M_{x,y} \tilde\Lap (\gh_n - \gh_m)(y) \\
      &= \sum_{x} \tilde \Lap (\gh_n - \gh_m)(x) (\gh_n - \gh_m)(x)    
      \qq\text{by \eqref{eqn:LapM=I}}.
    \end{align*}
    This is the inner product of two Cauchy sequences tending to 0 (by choice of \gh), and hence tends to 0.
    Since $\Lap_2$ and \sM are both self-adjoint, (ii) now follows from the spectral theorem. 
    %let $\gx \in \dom \sM$ and define $\gh := \sM\gx$, i.e., suppose that $(\gx_n)_{n=1}^\iy \ci \spn\{\gd_x\}$ and that there exists $\gx,\gh \in \ell^2(G)$ satisfying 
    %\linenopax
    %\begin{align*}%\label{eqn:}
    %  \lim_{k \to \iy} \|\gx - \gx_k\|_2 = 0
    %  \qq\text{and}\qq \lim_{k \to \iy} \|\gh - \sM \gx_k\|_2 = 0.
    %\end{align*}
    %\marginpar{\qq$\leftarrow$why?}
    %It follows that $\gh \in \dom \Lap_2$ and $\Lap_2 \gh = \gx$. \\
    %Hence, the operator $M$ on $\ell^2(G)$ satisfies $M\gx = \gh$ and $\Lap_2 M\gx = \gx$. It is immediate that $\Lap_2$ maps $\dom \Lap_2$ into $\dom M$ (again because $\Lap_2$ is banded), whence one also has $M \Lap_2 \gh = \gh$.
  \end{proof}
\end{lemma}

\begin{remark}%[Locally finite network/banded Laplacian]
\label{rem:banded}
In general, it is difficult to determine spectral properties of operators in Hilbert space from a representation of the operators in the form of an infinite matrix. This was noted by von Neumann in \cite{vN43,vN51}. However, restriction to the class of \emph{banded operators} allows one to obtain many explicit results; see \cite{Jor78}, for example. An infinite matrix is \emph{banded} iff each row and each column has only finitely many nonzero entries. %Banded matrices occur frequently in mathematical physics, and include Jacobi matrices, Heisenberg's position and momentum matrics $P,Q$, and polynomials in $P$ and $Q$. In the context of graphs and networks, banded matrices arise as graph Laplacians (or incidence matrices) for examples where each vertex has only finitely many neighbours.
  In the present context, the assumption of local finiteness of the network is equivalent to the bandedness of the Laplacian (on \HE or on $\ell^2(G)$); % is a banded operator (see Definition~\ref{def:banded-operator} and Definition~\ref{def:banded-operator} just below) because Definition~\ref{def:network} requires that $(\Graph,\cond)$ has finite degree. This condition is used only in the proof Lemma~\ref{thm:intertwined-powers}, and it is possible that it may be unnecessary.
  this hypothesis is used only for Lemma~\ref{thm:intertwined-powers}. It is quite possible that there may exist an alternative proof, in which case this hypothesis may turn out to be unnecessary.    
\end{remark}

\begin{defn}\label{def:Phi-ext}
  Define 
  \linenopax
  \begin{align}\label{eqn:Phi-ext}
    \widetilde \gF(\gx) := \lim_{n \to \iy} \gF(\gx_n),
    \qq \text{for any } \gx \in \dom q_M,
    %\widetilde \gF := \dom \Lap_2^{1/2} \cap \dom \Lap_2^{-1/2}, 
  \end{align}
  where $\left(\gx_n\right)_{n=1}^\iy \ci \dom\gF$ is any sequence for which $\lim_{n \to \iy} q_M(\gx_n-\gx) = 0$.
\end{defn}

\begin{lemma}\label{thm:Phi-tilde}
  The operator $\widetilde \gF$ is well-defined.  
  \begin{proof}
    Let $\gx \in \dom q_M$, and let $(\gx_n)_{n=1}^\iy$ be a sequence of finitely supported functions for which $\lim_{n \to \iy}q_M(\gx_n-\gx) = 0$.
    %tend to \gx in $\ell^2(G)$. 
    Then Lemma~\ref{thm:M-intertwines} gives
    \linenopax
    \begin{align*}%\label{eqn:}
      \|\gF(\gx_n - \gx_m)\|_\energy^2 
      &= \la\gx_n - \gx_m,M(\gx_n - \gx_m)\ra_2
      = q_M(\gx_n - \gx_m),
    \end{align*}
    which converges because $q_M$ is closed.  
  \end{proof}
\end{lemma}

\begin{remark}\label{rem:isometric-extension}
  Note that \gF is an isometry from $\spn\{\gd_x\}$ (equipped with the $q_M$-norm) into \HE; Lemma~\ref{thm:Phi-tilde} just emphasizes that this isometry is maintained under completion.
\end{remark}

\begin{comment}\label{thm:dom(qM)}
  The domain of $q_M$ is $\dom q_M = \{\gx \in \ell^2(G) \suth \|\widetilde\gF(\gx)\|_\energy<\iy\}.$ 
  \begin{proof}
    Combine Lemma~\ref{thm:M-intertwines} with Definition~\ref{def:Phi-ext}.
  \end{proof}
\end{comment}

\begin{defn}\label{def:q}
  For $u \in \dom \LapE$, define the quadratic form $r(u) := \la u, \Lap u\ra_\energy + \|u\|_\energy^2$,
  %on the domain %for $u \in \dom \LapE$ for which $\Lap u \in \HE$. That is, 
  %\linenopax
  %\begin{align}\label{eqn:dom-q}
  %  \dom q_\Lap := \{u \in \dom \LapE \suth \Lap u \in \HE\}.
  %\end{align}
  and denote the closure of this form (and its domain) by $q_\Lap$.
  %$q_\Lap$-completion of $\dom \LapE$ by
  %\linenopax
  %\begin{align}
  %  q_\Lap(u) := \la u, \Lap u\ra_\energy + \|u\|_\energy^2,
  %  \qq u \in \dom q_\Lap.
    %\{\lim_{n \to \iy} u_n \suth u_n \in \dom \LapE 
    %\;\text{ and } \lim_{m,n \to \iy} q_\Lap(u_m-u_n) = 0\},
  %\end{align}
  %where the first limit is taken with respect to the topology induced by \eqref{eqn:q}.
\end{defn}

\begin{defn}[Friedrichs extension]\label{def:Friedrichs}
  The \emph{Friedrichs extension} \LapF is the unique self-adjoint and non-negative   operator (with greatest lower bound 0) associated to $q_\Lap$ by Kato's theorem. 
\end{defn}

\begin{remark}\label{rem:Kato-qLap}
  Kato's theorem (Theorem~\ref{thm:Kato-rep}) gives $\dom \LapF^{1/2} = \dom q_\Lap$ and
  \begin{align}\label{eqn:D}
    q_\Lap(u) = \|\LapF^{1/2}u \|_\energy^2 + \|u\|_\energy^2,
    \qq\text{for } u \in \dom q_\Lap.
  \end{align}
\end{remark}

\begin{lemma}\label{thm:Fried-domain}
  The domain of the Friedrichs extension \LapF is  
  \linenopax
  \begin{align}\label{eqn:domLF}
    \dom \LapF = (\dom \LapE^\ad ) \cap (\dom q_\Lap).
  \end{align}
  \begin{proof}
    This follows from \cite[IV.3]{Kat95} or \cite[\S{}XII.5]{DuSc88}. Recall that convergence in energy implies pointwise convergence. 
  \end{proof}
\end{lemma}

\begin{defn}\label{def:Lap(-1/2)}
  %Let $\gx \in \dom \Lap_2$. 
  Since $\Lap_2$ is self-adjoint (see Remark~\ref{rem:essentially-self-adjoint}), we %use spectral theory to 
  define the operator
  %$\Lap_2^{-1/2}$ by : %let $\gm_\gx^{\ell^2}$ denote the spectral measure in the spectral resolution of $\Lap_2$. Then define
  \linenopax
  \begin{align}\label{eqn:Lap(-1/2)}
    \Lap_2^{-1/2} := \gp^{-1/2} \int_0^\iy t^{-1/2} e^{-t \Lap_2} dt
  \end{align}
  for those functions $\gx$ lying in the domain
  \linenopax
  \begin{align}\label{eqn:domLap(-1/2)}
    %\dom \Lap_2^{-1/2} := \{\gx \suth \gl^{-1/2} \in L^2(\bR_+,\gm_\gx^{\ell^2})\} 
    \dom \Lap_2^{-1/2} := \left\{\gx \suth \left(\gp^{-1/2} \int_0^\iy t^{-1/2} e^{-t \Lap_2} dt\right)(\gx) \in \ell^2(G)\right\}.
  \end{align}
  Observe that this integral converges because $0$ is not an eigenvalue of $\Lap_2$; recall that we consider only infinite networks. %We now define $\Lap_2^{-1/2}$ via functional calculus, by taking the square root of \eqref{eqn:Lap(-1)}. 
  See \cite[Ch.~XII]{DuSc88} or \cite[Ch.6--7]{Nel69}.
\end{defn}
  
The characterization of the Friedrichs domain extension given in Theorem~\ref{thm:Friedrichs} will require the following two lemmas.

\begin{lemma}\label{thm:D}
  $\dom \LapF^{1/2} = \widetilde\gF (\dom \Lap_2^{-1/2})$.
  \begin{proof}
    %\dom q_\Lap consists of limits of sequences $(u_n)_{n=1}^\iy$ which satisfy 
    In light of Remark~\ref{rem:Kato-qLap}, it suffices to show \smash{$\dom q_\Lap = \widetilde\gF (\dom \Lap_2^{-1/2})$}. 
    For any $u \in \dom q_\Lap$, one can find $(u_n)_{n=1}^\iy$ with $\lim q_\Lap(u_n-u_m)=0$ and $u_n = \gF\gx_n$ for $\gx_n \in \dom\gF$. Then
    \linenopax
    \begin{align}%
      q_\Lap(u_n-u_m)
      %&= \la u_n-u_m, \Lap(u_n-u_m)\ra_\energy
      %   + \|u_n-u_m\|_\energy^2 \\
      &= \la \gF(\gx_n - \gx_m), \Lap \gF(\gx_n - \gx_m)\ra_\energy
         + \|\gF(\gx_n - \gx_m)\|_\energy^2 \notag \\
      &= \|\gx_n - \gx_m\|_2 
        + \la\gx_n - \gx_m,M(\gx_n - \gx_m)\ra_2^2,
        \label{eqn:q_Lap-Cauchy}
    \end{align}
    by Lemma~\ref{thm:Phi-intertwines} and Lemma~\ref{thm:M-intertwines}.
    Now by Lemma~\ref{thm:Green} and Lemma~\ref{thm:Phi-tilde}, the convergence of \eqref{eqn:q_Lap-Cauchy} is equivalent to both $(\gx_n)_{n=1}^\iy$ and $(\Lap_2^{-1/2}\gx_n)_{n=1}^\iy$ being Cauchy in $\ell^2(G)$, but this means precisely that $\gx := \lim \gx_n \in \dom \Lap_2^{-1/2}$. 
    
    Conversely, if $\gx$ is the limit in $\ell^2(G)$ of a sequence $(\gx_n)_{n=1}^\iy \ci \dom\gF$, then observe that for $u_n = \gF\gx_n$, the same identity follows by Definition~\ref{def:q}.
  \end{proof}
\end{lemma}

\begin{comment}\label{rem:3-properties-of-completion}
  From the proof of Lemma~\ref{thm:D}, one can also see that $(u_n)_{n=1}^\iy$ is Cauchy in $q_\Lap$ if and only if the following three properties hold:
\linenopax
\begin{align}
  u_n = \gF(\gx^{(n)}), &\qq\text{ for some } \gx^{(n)} \in \dom\gF, \label{eqn:un}\\
  \|u_m-u_n\|_\energy &\limas{m,n} 0, \qq\text{and} \label{eqn:unE}\\
  \|\gx^{(m)}-\gx^{(n)}\|_2 &\limas{m,n} 0. \label{eqn:gx2}
\end{align} 
\end{comment}

\begin{lemma}\label{thm:adjoint-domains-identity}
  $\dom \LapE^\ad = \widetilde\gF(\dom \Lap_2^{1/2})$.
  \begin{proof}
    To begin, we show that any element $u \in \dom \LapE^\ad$ can be written as $u = \widetilde\gF(\gx)$ for some $\gx \in \dom \widetilde \gF$.
    By Definition~\ref{def:domLapE} and Definition~\ref{def:Phi-ext}, $u \in \dom \LapE^\ad$ iff there exists a $C < \iy$ for which 
    \linenopax
    \begin{align}\label{eqn:C-bound1}
      \left|\la\LapE \widetilde\gF\gh, u\ra_\energy\right|^2 
      \leq C \|\widetilde\gF\gh\|_\energy^2, 
      \qq\text{for all } \gh \in \dom \widetilde\gF.
    \end{align}
    Combining Lemma~\ref{thm:M-intertwines}, Lemma~\ref{thm:Phi-tilde}, and Definition~\ref{def:qM} gives $\|\widetilde\gF\gh\|_\energy^2 = \|\Lap_2^{-1/2}\gh\|_2^2$, so that \eqref{eqn:C-bound1} is equivalent to
    \linenopax
    \begin{align}\label{eqn:C-bound2}
      \left|\la \LapE \widetilde\gF\gh, u\ra_\energy\right|^2 
      \leq C \| \Lap_2^{-1/2} \gh\|_2^2, 
      \qq\text{for all } \gh \in \dom \widetilde\gF.
    \end{align}
    Now Riesz duality gives a $\gx \in \ell^2(G)$ for which $\la \Lap_2^{-1/2} \gh, \gx\ra_{\ell^2} = \la \LapE \widetilde\gF\gh, u\ra_\energy$, so with
    %\eqref{eqn:C-bound2} is equivalent to
    %\linenopax
    %\begin{align}\label{eqn:C-bound3}
    %  \left|\la \gh, \gx\ra_{\ell^2}\right|^2 
    %  \leq C \| \Lap_2^{-1/2} \gh\|_2^2, 
    %  \qq\text{for all } \gh \in \dom \widetilde\gF.
    %\end{align}
    %By Lemma~\ref{thm:Phi-intertwines}, this is equivalent to 
    %\linenopax
    %\begin{align}\label{eqn:C-bound2}
    %  \left|\la\gh,\gx\ra_2\right|^2 
    %  \leq C \left|\la\gh,M\gh\ra_2\right|^2
    %  = C \|\Lap_2^{-1/2}\gh\|_2^2, 
    %  \qq\text{for all } \gh \in \dom q_M,
    %\end{align}
    %where $\Lap_2^{-1}$ is defined in \eqref{eqn:Lap(-1/2)}.
    %Setting 
    $\gr = \Lap_2^{-1/2} \gh$, we can rewrite \eqref{eqn:C-bound2} as
    \linenopax
    \begin{align*}%\label{eqn:}
      \left|\la \gr, \gx\ra_2\right|^2 
      \leq C \|\gr\|_2^2, 
      \qq\text{for all } \gr \in \Lap_2^{-1/2}\left(\dom q_M\right),
    \end{align*}
    which means $\gx \in \dom (\Lap_2^{1/2})^\ad = \dom \Lap_2^{1/2}$, since $\Lap_2$ is self-adjoint by \cite[Lem.~2.7 and Thm.~2.8]{SRAMO}; see Remark~\ref{rem:essentially-self-adjoint}.
  \end{proof}
\end{lemma}

\begin{theorem}[Friedrichs characterisation]
  \label{thm:Friedrichs}
  The Friedrichs extension is given by 
  \linenopax
  \begin{align}\label{eqn:domLF-characterization}
    \dom \LapF = \widetilde\gF(\dom \Lap_2^{1/2}\cap \dom \Lap_2^{-1/2}).
  \end{align}
  \begin{proof}
    Starting with \eqref{eqn:domLF}, applying Lemma~\ref{thm:D} and Lemma~\ref{thm:adjoint-domains-identity} gives
    \linenopax
    \begin{align}
      \dom \LapF 
      = \dom \LapE^\ad \cap \dom q_\Lap
      &= \widetilde\gF(\dom \Lap_2^{1/2}) \cap \widetilde\gF(\dom \Lap_2^{-1/2})         
      \label{eqn:thingyquestion1}
    \end{align}
    Suppose that $u \in \dom \LapF$ can be written as $u=\gF(\gx)$ for $\gx \in \dom \Lap_2^{1/2}$ and as $u=\gF(\gh)$ for $\gh \in \dom \Lap_2^{-1/2}$ . Then
    \linenopax
    \begin{align*}%\label{eqn:rep-equality}
      \la \gd_x, \gF\gx\ra_\energy 
      = \sum_{y \in G} \gx_y \la \gd_x, w_y \ra_\energy
      = \sum_{y \in G} \gh_y \la \gd_x, w_y \ra_\energy
      = \la \gd_x, \gF\gh\ra_\energy
    \end{align*}
    Therefore, $\gx = \gh$ and $\widetilde\gF$ preserves intersections, so \eqref{eqn:thingyquestion1} is equal to \eqref{eqn:domLF-characterization}.
    %Suppose $u \in \dom \LapF$. Then $u \in \dom \LapE^\ad$ by \eqref{eqn:domLF}, so Lemma~\ref{thm:adjoint-domains-identity} implies $\gx \in \dom \Lap_2^{1/2}$. Now suppose $u \in \widetilde\gF(\dom \Lap_2^{1/2}) \cap \HE$, so $u = \gF\gx$ for some $\gx \in \dom \Lap_2^{1/2}$ and Lemma~\ref{thm:adjoint-domains-identity} implies $u \in \dom \LapE^\ad$. Finally, note that if $u \in \HE$, then there is a representation of $u$ as $u = \sum_{x \in X} \gx_x w_x$. If $(F_n)_{n=1}^\iy$ is an exhaustion of $X$, then 
    %\linenopax
    %\begin{align*}%\label{eqn:}
    %  u_n := \sum_{x \in F_n} \gx_x w_x
    %\end{align*}
    %converges to $u$ in \HE, as $n \to \iy$. 
    %It follows from Lemma~\ref{thm:D} that $u \in \dom q_\Lap$, and hence in $\dom \LapF$.
  \end{proof}
\end{theorem}

%\begin{cor}\label{thm:F-domain-is-Dirichlet}
%  If $u \in \dom \LapF$, then $u$ vanishes on $\bd G$. 
%  \begin{proof}
%    \marginpar{Proof?}
%    Suppose $u = \gF \gx \in \dom \LapF$. Then ...
%  \end{proof}
%\end{cor}

%%!TEX root = Friedrichs.tex

\section{Relating the Friedrichs extension \LapF on \HE to $\Lap_2$}
\label{sec:Fried-history}

The main result of this section is Theorem~\ref{thm:spectral-resolutions}, in which we show that the spectral measures of \LapF and $\Lap_2$ are mutually absolutely continuous in the complement of $\gl=0$ and compute the Radon-Nikodym derivative. 

\begin{defn}\label{def:spectral-measures}
  For $u \in \dom \LapF$, let $\gm_{u}^{\sF}$ denote the spectral measure in the spectral resolution of \LapF, and for $\gx \in \dom \Lap_2$, again let $\gm_\gx^{\ell^2}$ denote the spectral measure in the spectral resolution of $\Lap_2$.
\end{defn}

\begin{theorem}\label{thm:spectral-resolutions}
  For $\gx \in \dom \gF$, the spectral measures of \LapF and $\Lap_2$   
  %Definition~\ref{def:spectral-measures} 
  are related by
  \linenopax
  \begin{align}\label{eqn:spectral-resolutions}
    \gl d\gm_{\gF\gx}^{\sF} = d\gm_\gx^{\ell^2},
  \end{align}
  where \gl is the spectral parameter.
  In particular, $d\gm_{\gF\gx}^{\sF}$ and $d\gm_\gx^{\ell^2}$ are mutually absolutely continuous on $(0,\iy)$ with Radon-Nikodym derivative $\gl$.
  \begin{proof}
    For $\gx \in \dom\gF$, 
    \linenopax
    \begin{align*}%\label{eqn:}
      \la \gF\gx, \Lap^{n+1}\gF\gx\ra_\energy
      = \la \gF\gx, \Lap\gF\Lap^{n}\gx\ra_\energy
      = \la \gx, \Lap^{n}\gx\ra_2
    \end{align*}
    follows by applying Lemma~\ref{thm:intertwined-powers} and then Lemma~\ref{thm:Phi-intertwines} (with $\gh = \Lap^n \gx$). Note that \gx has finite support, and therefore so does $\Lap^{n}\gx$ for any $n$; see the proof of Lemma~\ref{thm:Green}. This identity also holds for operators with larger domains, so
    \linenopax
    \begin{align}%\label{eqn:}
      \la \gF\gx, \LapF^{n+1}\gF\gx\ra_\energy
      = \la \gx, \Lap_2^{n}\gx\ra_2,
      \qq\text{for all } \gx \in \dom \Lap_2^n.
    \end{align}
    If $P_2$ denotes the projection-valued measure in the spectral resolution of $\Lap_2$ and $P_\sF$ denotes the projection-valued measure in the spectral resolution of \LapF, then the spectral theorem gives
    %\linenopax
    %\begin{align}%\label{eqn:}
    %  \left\|\int_0^\iy \gl^{n+1} P_\sF(d\gl)\gF\gx\right\|_\energy^2 
    %  = \left\|\int_0^\iy \gl^{n} P_2(d\gl)\gx \right\|_2^2, 
    %\end{align}
    %which by the orthogonality of $P_\sF$ and $P_2$ can be rewritten
    \linenopax
    \begin{align}%\label{eqn:}
      \int_0^\iy \gl^{n+1} \|P_\sF(d\gl)\gF\gx\|_\energy^2 
      = \int_0^\iy \gl^{n} \|P_2(d\gl)\gx\|_2^2. 
    \end{align}
    Considering the above as an identity for monomials $\gl^n$, it is clear the measures %a standard application of measure theory to see that 
    \linenopax
    \begin{align*}%\label{eqn:}
      \gl d\gm_{\gF\gx}^{\sF}(\gl) := \gl \|P_\sF(d\gl)\gF\gx\|_\energy^2
      \qq\text{and}\qq
      d\gm_\gx^{\ell^2}(\gl) := \|P_2(d\gl)\gx\|_2^2
    \end{align*}
    have the same moments. Since $\Lap_2$ is essentially self-adjoint, the corresponding moment problem is determinate, i.e., these moments determine the measures uniquely (see, e.g., \cite{Akhiezer, AkhiezerGlazman, Fuglede83}). %By RieszÕs theorem (see, e.g., \cite[Thm.~2.14]{Rudin}) we can identify measures with linear functionals on continuous functions. So the aim to get from monomials \gd^n to all Borel functions. Or, by Riesz, equivalently, to all continuous functions suitably restricted. This is the key step in the proof of eq (5.1), so the theorem itself. The thing is that the moments do not always determine measures uniquely. However in this case, we are ok on account of our l^2-theorem: Lap in l^2  is essentially selfadjoint.
    Consequently, we have 
    \linenopax
    \begin{align*}%\label{eqn:}
      \int_0^\iy \gl f(\gl) \,d\gm_{\gF\gx}^{\sF} (\gl)
      = \int_0^\iy f(\gl) \, d\gm_\gx^{\ell^2}(\gl)
    \end{align*}
    for any bounded Borel function $f$. (Note that $f$ is required to be bounded because the operators $\Lap_2, \LapF$ may not be.) This completes the proof of \eqref{eqn:spectral-resolutions}.
  \end{proof}
\end{theorem}

\begin{remark}\label{rem:spectral-resolution-consequences}
  Note that we are only considering infinite networks here, so $0$ is never an eigenvalue of $\Lap_2$. Thus, \eqref{eqn:spectral-resolutions} states that the spectra of \smash{$\gm_{\gF\gx}^{\sF}$} and \smash{$\gm_\gx^{\ell^2}$} must agree on the support of the Radon-Nikodym derivative, i.e., up to an eigenvalue at $0$ for  $\gm_{\gF\gx}^{\sF}$ (which is present precisely in the case $\Harm \neq 0$).
  Several useful facts follow immediately from $\spec_{\HE} \LapF \less \{0\} = \spec_{\ell^2} \Lap_2$ and \eqref{eqn:spectral-resolutions}. For example,  
  %the formulas
  %\linenopax
  %\begin{align}\label{eqn:spectral-measures-normsquares}
  %  \gl \|P_\sF(d\gl) \gF \gx\|_\energy^2 = \|P_2(d\gl)\gx\|_2^2,
  %  \qq\text{for all } \gx \in \dom\gF,
  %\end{align}
  %and
  \linenopax
  \begin{align}\label{eqn:spectral-measures-Phi(gx)}
    \|\gF\gx\|_\energy^2  = \| \Lap_2^{-1/2} \gx\|_2^2,
    \qq\text{for all } \gx \in \dom \Lap_2^{-1/2},
  \end{align} 
  %Formula \eqref{eqn:spectral-measures-Phi(gx)} 
  which allows us to prove transience of the integer lattice networks in Theorem~\ref{thm:monopoles-on-Zd}.
  Corollary~\ref{thm:bounded} is another useful consequence of this fact.
  
  It also follows that for transient networks, the Dirac measure at $\lambda = 0$ contributes to the spectral resolution of the Friedrichs extension of $\Delta_{\mathcal{E}}$ but not to that of the self-adjoint $\ell^2$ Laplacian.
\end{remark}

\begin{cor}\label{thm:bounded}
  \LapE is bounded if and only if $\Lap_2$ is bounded.
  \begin{proof}
    %By Theorem~\ref{thm:spectral-resolutions}, 
    Note that \LapF is bounded if and only if $\|\LapF\| = \sup \spec \LapF < \iy$. It is immediate from \eqref{eqn:spectral-resolutions} that $\sup \spec \LapF = \sup \spec \Lap_2$, and hence  $\|\LapF\| < \iy$ is equivalent to $\|\Lap_2\| = \sup \spec \Lap_2 < \iy$.    
    Since \LapF is an extension of \LapE (which coincides with \LapE whenever \LapE is bounded), the result is immediate.
  \end{proof}
\end{cor}

\begin{cor}\label{thm:Phi(gx)-in-HE}
  Let $\gx \in \ell^2(G)$. Then $\gF(\gx) \in \HE$ if and only if $\gx \in \ran \Lap_2^{1/2}$.
  \begin{proof}
    By the spectral theorem,
    \linenopax
    \begin{align}\label{eqn:Phi(gx)-in-HE-der1}
      \|\gF\gx\|_\energy^2
      &= \int_0^\iy \| P_\sF(d\gl) \gF\gx\|_\energy^2
       = \int_0^\iy \frac1\gl \| P_2(d\gl)\gx\|_2^2
       = \| \Lap_2^{-1/2} \gx\|_2^2,
    \end{align}
    where the second equality comes by \eqref{eqn:spectral-resolutions}, since $\frac1\gl = \left(\gl^{-1/2}\right)^2$. Then \eqref{eqn:Phi(gx)-in-HE-der1} is finite if and only if $\gh = \Lap_2^{-1/2} \gx \in \ell^2(G)$, that is, %$\Lap_2^{-1/2} \gx = \gh$, i.e., 
    $\gx \in \ran \Lap_2^{1/2}$.
  \end{proof}
\end{cor}

%\input{conductances}
%\input{currents}
%%!TEX root = Friedrichs.tex

\section{Applications to effective resistance}

The main result of this section is Theorem~\ref{thm:R(x,y)}, a corollary to Theorem~\ref{thm:spectral-resolutions} which may be compared with Remark~\ref{rem:wired-resistance}. Recall the monopole notation \monov from Remark~\ref{thm:minimal-monopole}. We use \Pfin to denote the orthogonal projection of $v_x$ to \Fin with respect to $\la \cdot, \cdot \ra_\energy$; see Theorem~\ref{thm:HE=Fin+Harm}. Thus  $f_x := \Pfin v_x$ and $\monof := f_x + w_o = \Pfin \monov$.

\begin{defn}\label{def:R(x)}
  Denote the \emph{free effective resistance} between $x$ and $y$ by 
  \linenopax
  \begin{align}\label{eqn:R(x)}
    R^F(x,y) := \energy(\monov - \monoy) = \energy(v_x - v_y),
  \end{align}
  for \monov as in \eqref{eqn:w_x}.  
  Denote the \emph{wired effective resistance} between $x$ and $y$ by 
  \linenopax
  \begin{align}\label{eqn:RW(x)}
    R^W(x,y) := \energy(f_x - f_y) = \energy(\monof - \monof[y]).
  \end{align}  
\end{defn}

\begin{remark}\label{rem:wired-resistance}\label{rem:RF-vs-RW}
  Several alternative and equivalent formulations of the free and wired resistances are collected in \cite{ERM}. 
  It turns out that $R^F$ and $R^W$ are metrics on $G$; for details, see \cite{ERM,Lyons,Kig03}. Note that $R^F(x,y) \geq R^W(x,y)$ in general, and that strict inequality holds if and only if $\Harm \neq 0$. %, since $v_x=f_x$ (and $\monov = \monof$) when $\Harm=0$. 
\end{remark}

The following lemma combines results from \cite[Lem.~6.9]{bdG} and \cite[Lem.~2.22]{DGG} and will be useful in the sequel.

\begin{lemma}
  \label{thm:energy-kernel-is-real-and-bounded}
  Every $\monof$ is \bR-valued, with $w_x(y) - w_x(o) >0$ for all $y \neq o$. Moreover, every $w_x$ is bounded, with $\|w_x\|_\iy \leq R^F(x,o)$ (see \eqref{eqn:R(x)}). % and $\|\monov\|_\iy = $.
\end{lemma}

\begin{defn}\label{def:p(x,y)}
  The probabilities $p(x,y) := {\cond_{xy}}/{\cond(x)}$ define a random walk $(X_n)_{n=0}^\iy$ on the network by $\prob[X_{n+1}=y|X_n=x]=p(x,y)$. Here $X_n$ is a $G$-valued random variable giving the location of the random walker at time $n$. Then let
  \linenopax
  \begin{align}\label{eqn:Prob[x->y]}
    \prob[x \to y] := \prob_x(\gt_y < \gt_x^+) 
  \end{align}
  be the probability that the random walk started at $x$ reaches $y$ before returning to $x$. In \eqref{eqn:Prob[x->y]}, $\gt_z := \min\{n \geq 0 \suth X_n = z\}$ is the hitting time of $z$ and $\gt_z^+ := \max\{\gt_z,1\}$.
\end{defn}

\begin{cor}[{\cite[Cor.~3.13 and Cor.~3.15]{ERM}}]
  \label{thm:c(x)R(x)=Prob[x->o]}
  For any $x \neq y$, one has
  \linenopax
  \begin{equation}\label{eqn:c(x)R(x)=Prob[x->o]}
    R^F(x,o) = \frac1{\cond(o) \prob[o \to x]}.
  \end{equation}
\end{cor}

\begin{lemma}\label{thm:lower-bound-on-R}
  If \LapE is bounded on \HE, then $(G,R^F)$ is uniformly discrete, i.e., 
  \linenopax
  \begin{align}\label{eqn:lower-R-bound}
    R^F(x,y) \geq \frac{2}{\|\LapE\|}.
  \end{align}
  \begin{proof}
    Since the inequality $\la u, \Lap u \ra_\energy \leq \|\LapE\| \cdot \|u\|_\energy^2$ holds for all $u \in \dom \LapE = \HE$, apply it to $u = \monov - \monoy$ (with $x \neq y$) to obtain
    \linenopax
    \begin{align}\label{eqn:Steinhaus-bound-applied}
      \la \monov - \monoy, \Lap(\monov - \monoy) \ra_\energy 
      \leq \|\LapE\| \cdot \|\monov - \monoy \|_\energy^2 .
    \end{align}
    Note also that Lemma~\ref{thm:<delta_x,v>=Lapv(x)} gives
    \linenopax
    \begin{align}
      \la \monov - \monoy, \Lap(\monov - \monoy) \ra_\energy 
      &= \la \monov - \monoy, \gd_x - \gd_y \ra_\energy \notag \\
      &= \la \monov, \gd_x\ra_\energy - \la \monov, \gd_y\ra_\energy - \la \monoy, \gd_x\ra_\energy + \la \monoy, \gd_y\ra_\energy \notag \\
      &= \Lap \monov(x) - \Lap \monov(y) - \Lap \monoy(x) + \Lap \monoy(y) \notag \\
      &= 1 - 0 - 0 + 1 = 2.
      \label{eqn:Lapcalc}
    \end{align}
    Combining \eqref{eqn:Steinhaus-bound-applied}, \eqref{eqn:Lapcalc}, and \eqref{eqn:R(x)} gives \eqref{eqn:lower-R-bound}.
  \end{proof}
\end{lemma}

\begin{remark}\label{rem:0}
  Note that $0$ is never an eigenvalue of $\Lap_2$.%
    \footnote{It is well-known that; the only harmonic function in $\ell^2(G)$ on an infinite network $G$ is the constant function $0$; %. For example, it follows from \cite[Lem.~5.5]{DGG} that any harmonic function $h$ on an infinite network is bounded away from 0 on an infinite set, and hence $h$ cannot lie in $\ell^2(G)$. See also 
    see \cite{Lyons,Woess09,Soardi94,DGG} and elsewhere.} 
It follows that we can only apply \eqref{eqn:spectral-resolutions} in the orthogonal complement of \Harm (since the formula may not hold at $\gl=0$), so there is no analogue of Theorem~\ref{thm:R(x,y)} for $R^F$.
\end{remark}

\begin{theorem}\label{thm:R(x,y)}
  For an infinite network $G$, the wired effective resistance $R^W(x,y)$ is %given by
  \linenopax
  \begin{align}\label{eqn:R(x,y)}
    R^W(x,y) = \int_0^\iy \frac1\gl \|P_2(d\gl)(\gd_x - \gd_y)\|_2^2
  \end{align}
  \begin{proof}
    $R^W(x,y) = \| \monof -  \monof[y]\|_\energy^2$, we have
    \linenopax
    \begin{align*}%\label{eqn:R(x,y)-der1}
      R^W(x,y)
      %= \| \monov -  \monov[y]\|_\energy^2 
      = \int_0^\iy \| P_\sF(d\gl) (\monof -  \monof[y])\|_\energy^2 
      &= \int_0^\iy \| P_\sF(d\gl) \Pfin \gF(\gd_x-\gd_y)\|_\energy^2. 
    \end{align*}
    Note that \Pfin is the projection to the orthocomplement of $\ker \LapF$, so we can remove it\footnote{Let $T = T^\ad$ be a self-adjoint operator densely defined on the Hilbert space $\sH = \sH_0 \oplus \sH_1$, where $\sH_0 = \ker T$. If $P_1$ is projection to $\sH_1$, then $TP_1=T$.}:
    \linenopax
    \begin{align*}%\label{eqn:R(x,y)-der1}
      R^W(x,y)
      &= \int_0^\iy \| P_\sF(d\gl) \gF(\gd_x-\gd_y)\|_\energy^2
      = \int_0^\iy \frac1\gl \|P_2(d\gl)(\gd_x - \gd_y)\|_2^2,
    \end{align*}
    where the last equality follows by Theorem~\ref{thm:spectral-resolutions}.
  \end{proof}
\end{theorem}

\begin{cor}\label{thm:R(x,y)-spectral-gap}
  On an infinite network $G$, If $\Lap_2$ has a spectral gap, then the wired effective resistance is bounded. More precisely, if $\spec \Lap_2 \ci [\gg,\iy)$, then
  \linenopax
  \begin{align}\label{eqn:R(x,y)-spectral-gap}
    R^W(x,y) \leq \frac2\gg \;,
    \qq\text{for all } x,y \in G.
  \end{align}
  \begin{proof}
    In this case, \eqref{eqn:R(x,y)} gives
    \linenopax
    \begin{align*}%\label{eqn:R(x,y)-spectral-gap-der1}
      R^W(x,y)
      &= \int_\gg^\iy \frac1\gl \|P_2(d\gl)(\gd_x - \gd_y)\|_2^2 
      \leq \frac1\gg \int_\gg^\iy \|P_2(d\gl)(\gd_x - \gd_y)\|_2^2
      = \frac1\gg \|\gd_x - \gd_y\|_2^2,
    \end{align*}
    and $\|\gd_x - \gd_y\|_2^2 \leq 2$ with strict inequality iff $x=y$, in which case $R(x,y)=0$.
  \end{proof}
\end{cor}

%\begin{conj}\label{rem:boundedness-conjecture}
%  We expect that the converse of Corollary~\ref{thm:R(x,y)-spectral-gap} holds as well, i.e., that a bound on $R^W$ occurs only in the presence of a spectral gap for $\Lap_2$.
%\end{conj}

%%!TEX root = Friedrichs.tex

\section{Examples}
\label{sec:examples}

In this section, we apply our results to obtain concrete formulas for some common and well-studied examples, including trees and integer lattices. %The underlying group structure of this example permits the use of Fourier theory to solve certain equations, and these can be readily translated into explicit spectral representations. %Throughout this section, we abuse notation and for $v \in \HE$, we write $v$ for the representative $v:G \to \bR$ which vanishes at the origin, i.e., $v(o)=0$. 

\subsection{The binary tree}
  \label{sec:tree2}
Consider the network $(\bT_2, \one)$, the binary tree with all edges having conductance 1. The vertices of this network can be labeled with finite words on the symbol set $\{0,1\}$, using $\es$ to denote the empty word, which we take as the origin, i.e. $o=\es$. For a vertex of the tree $x \in \{0,1\}^k$, let $|x|:=k$ (with $|o|=0$).
%\linenopax
%\begin{align}\label{eqn:binary-tree}
%  \xymatrix@R=-5pt@C=40pt{
%  \bT_2 & & \vertex{11} \ar@{:}[r] &\\
%  & \vertex{1} \ar@{-}[ur] \ar@{-}[r] & \vertex{10} \ar@{:}[r] &\\
%  \vertex{\es} \ar@{-}[ur] \ar@{-}[dr] & &  \\
%  & \vertex{0} \ar@{-}[dr] \ar@{-}[r] & \vertex{01} \ar@{:}[r] &\\
%  & & \vertex{00}\ar@{:}[r] &
%  } 
%\end{align}
Using symmetry and elementary calculations, it is easy to check that $w(x):=2^{-|x|}$ is the unique energy-minimizing monopole on $\bT_2$. 
On any tree with conductance function $c=\one$, it is straightforward to see that the free resistance $R^F(x,y)$ coincides with combinatorial (shortest-path) distance, and is hence unbounded. One way to see this is to show that $v_x-v_y$ has a representative $u$ defined as follows: let $\gg$ be the shortest path from $y$ to $x$, and for $s \in \gg$, define $u(s)$ to be the number of edges between $s$ and $y$. The, for $s \notin \gg$, let $s_\gg$ be the unique closest point of \gg to $s$, and define $u(s) = u(s_\gg)$. Now $u$ increases by 1 with each step along \gg and is locally constant outside of \gg. One can check that $u$ satisfies the reproducing property required by $v_x-v_y$, and it is immediate that $\energy(u)<\iy$. 
However, in \cite{DuJo10}, the spectral gap for $\Lap_2$ on this example is computed to be $3-2\sqrt2$, and so it follows from Corollary~\ref{thm:R(x,y)-spectral-gap} that 
\linenopax
\begin{align}\label{eqn:RW-bounded-in-tree}
  R(x,y) \leq \frac{2}{3-2\sqrt2},
  \qq \text{ for all } x,y \in G.
\end{align}
This result appears to be new to the literature. %See Figure~\ref{fig:tree-exhaustions} for intuition regarding \eqref{eqn:RW-bounded-in-tree}.

\subsection{Expansive networks for which $\Harm=0$}
Suppose we now consider the network $(\bT_2 \times \bZ,\one)$ formed by taking the Cartesian product of $(\bT_2, \one)$ with the 1-dimensional integer lattice $(\bZ,\one)$. 

\begin{comment}\label{def:cartesian-product}
  If $(X,a)$ and $(Y,b)$ are both networks, the \emph{Cartesian product} is defined to have vertex set $X \times Y$ and edges defined by 
  \linenopax
  \begin{align*}%\label{eqn:}
    c_{(x,y),(s,t)} 
    = \begin{cases}
      a_{x,s}\,, &\text{if } y=t, \\ 
      b_{y,t}\,, &\text{if } x=s, \\ 
      0, &\text{else}.
    \end{cases}
  \end{align*}
\end{comment}

%Thus, each vertex in $(\bT_2 \times \bZ,\one)$ may be written $(\gw,n)$, where $\gw \in \{0,1\}^k$ for some $k \in \bN$, and $n \in \bZ$. The vertex $(\gw,n)$ has three neighbours in $\bT_2 \times n$, one in $\bT_2 \times (n+1)$, and one in $\bT_2 \times (n-1)$. One can show that $\Harm=0$ for this network, and in fact for any similar network, via the following lemma.

\begin{comment}[{\cite[Ex.~9.7]{Lyons}}]
  \label{thm:products-have-Harm=0}
  $\Harm=0$ for any network which is the Cartesian product of two infinite locally finite networks with $c=\one$.
\end{comment}

\begin{defn}\label{def:strong-isoperimetric-inequality}
   The network $(G,c)$ satisfies a \emph{strong isoperimetric inequality} iff there exists $\gd > 0$ such that for any finite vertex subset $S$, one has
  \linenopax
  \begin{align}\label{eqn:strong-isoperimetric-inequality}
    \frac{|\del S|}{|S|} \geq \gd > 0,
    \qq \text{where } |S| = \sum_{x \in S} 1 \q\text{and}\q
    |\del S| = \sum_{(xy) \in \del S} c_{xy},
  \end{align}
  and $\del S$ is the set of edges with exactly one end (vertex) in $S$. The infimum of $\frac{|\del S|}{|S|}$ (taken over all nonempty finite subsets $S \ci G$) is called the \emph{expansion constant}. 
\end{defn}
  
It is well-known that $\Lap_2$ has a spectral gap (i.e., $\inf\{\spec \Lap_2\} > 0$ if and only if \eqref{eqn:strong-isoperimetric-inequality} is satisfied. Networks satisfying these equivalent properties are called \emph{expanders}.%
  \footnote{For more information on this equivalence and its connections to mixing times of Markov chains, Kazhdan's property $T$ (or more precisely, property \gt), and other fascinating topics, we refer the reader to \cite[\S6]{Chung} or \cite{Woess09,Lyons,Expanders}; see also the various works of Lubotzky, Z\.uk, Bourgain, Gamburd, and Sarnak.}
  To see that $\bT_2 \times \bZ$ satisfies \eqref{eqn:strong-isoperimetric-inequality}, first observe that any regular tree of degree $d > 2$ satisfies \eqref{eqn:strong-isoperimetric-inequality} with $\gd=d-1$. Next, use the fact that the expansian constant for a Cartesian product is the sum of the expansion constants for the two factor networks.
Furthermore, it follows from \cite[Ex.~9.7]{Lyons} that $\Harm=0$ for $(\bT_2 \times \bZ,\one)$, and so %$R^F(x,y) = R^W(x,y)$. %it is also not hard to see that the network $\bT_2 \times \bZ$ satisfies a strong isoperimetric inequality.
%The expansion property of $\bT_2 \times \bZ$ ensures a spectral gap $\gg>0$, and so 
one has $R^F(x,y) = R^W(x,y) \leq \frac2\gg$. %for this network, upon combining Lemma~\ref{thm:products-have-Harm=0} with Corollary~\ref{thm:R(x,y)-spectral-gap}. 
It is clear that these considerations hold more generally than in this example.

\begin{lemma}\label{thm:R(x,y)-bounded-on-expanders}
  Let $(G,c)$ be an infinite network which satisfies a strong isoperimetric inequality and for which $\Harm=0$. Then the (necessarily unique) effective resistance on $(G,c)$ is bounded.
\end{lemma}

For the example of the binary tree $(\bT_2,\one)$, note that $\Harm \neq 0$ and $R^F$ is unbounded, but also that $R^W$ is bounded as in \eqref{eqn:RW-bounded-in-tree}.

\begin{remark}[Gel'fand spaces, and the 1-point compactification of $(G,R^W)$]
  \label{rem:Gel'fand-spaces}
  Denote the collection of bounded functions of finite energy by
  \linenopax
  \begin{equation}
    \algE := \{u \in \HE \suth u \text{ is bounded}\}.
  \end{equation}
  Define multiplication on $\algE$ by the pointwise product $(u_1 u_2)(x) := u_1(x) u_2(x),$ and let the norm on $\algE$ be given by $\|u\|_\sA := \|u\|_\iy + \|u\|_\energy$. With these definitions, it is shown in \cite[Lem.~5.5]{Multipliers} that $(\algE,\| \cdot \|_\sA)$ is a Banach algebra.
  
  For a Banach algebra \sA, the associated \emph{Gel'fand space} is the spectrum $\spec(\sA)$ realized as either the collection of maximal ideals of \sA or as the collection of multiplicative linear functionals on \sA. See \cite{Arveson:spectral-theory,Arveson:invitation-to-Cstar}. In \cite[Thm.~5.12]{Multipliers}, it is shown that if $\Harm = 0$, then the 1-point compactification of $(G,R^W)$ coincides with the Gel'fand space of $\sA_\energy$.
\end{remark}

\subsection{The integer lattices \bZd}
  \label{sec:Zd}
  \noindent Consider the $d$-dimensional integer lattice network $(\bZd, \one)$ with vertices 
  \linenopax
  \begin{align}\label{eqn:}
    \bZ^d = \{x=(x_1,\dots,x_d) \suth x_i \in \bZ, i=1,\dots,d\}
  \end{align}
  and with unit-conductance edges between nearest neighbours, that is,
  \linenopax
  \begin{align}\label{eqn:Zd-c}
    c_{xy} = \begin{cases}
      1, & y=x+\ge_k \text{ for some } k=1,\dots,d,\\
      0,& \text{else,}
    \end{cases}
  \end{align} 
  where $\ge_k = [0,\dots,0,1,0,\dots,0]$ has the 1 in the \kth slot. Let $o=0=(0,\dots,0)$.
%\end{exm}

%\subsubsection{Fourier theory of \Lap on $(\bZd,\one)$}

The following result is well-known; see \cite{Soardi94}, for example.

\begin{lemma}\label{thm:Fourier-transform-of-Lap-on-Zd}
  On the network $(\bZd, \one)$, the Fourier transform of \Lap is multiplication by 
  \linenopax
  \begin{align}\label{eqn:S(t)}
    S(t) = S(t_1,\dots,t_d) = 4 \sum_{k=1}^d \sin^2\left(\tfrac{t_k}2\right).
  \end{align}
\end{lemma}

\begin{lemma}\label{thm:vx-on-Zd}
  Let $\{v_x\}_{x \in \bZ^d}$ be the energy kernel on the integer lattice $\bZd$ with $\cond=\one$. Then for $y \in \bZd$,
  \linenopax
  \begin{equation}\label{eqn:vx-on-Zd}
    v_x(y) = \frac1{(2\gp)^d} \int_{\bT^d} \frac{\cos((x-y)\cdot t) - \cos(y\cdot t)}{S(t)} \,dt,
  \end{equation}
  where $dt$ is Haar measure on the $d$-torus \bTd.
  \begin{proof}
    Under the Fourier transform, Lemma~\ref{thm:Fourier-transform-of-Lap-on-Zd} indicates that the equation $\Lap v_x = \gd_x - \gd_o$ becomes $S(t) \hat v_x = e^{\ii x\cdot t} - 1$, whence
    \linenopax
    \begin{align}\label{eqn:vx-on-Zd-as-exponential}
      v_x(y) = \frac1{(2\gp)^d} \int_{\bT^d} e^{-\ii y \cdot t} \frac{e^{\ii x\cdot t} - 1}{S(t)} \,dt.
    \end{align}
    Since $v_x$ is \bR-valued, the result follows.
  \end{proof}
\end{lemma}

\begin{remark}\label{rem:Harm-on-Zd}
  It is known that no nonconstant harmonic functions of finite energy exist on the integer lattices, and hence the free and wired resistance metrics on $(\bZd,\one)$ coincide \cite{Woess00, Woess09, Lyons, ERM}; see also Definition~\ref{def:R(x)} and Remark~\ref{rem:RF-vs-RW}. Hence, we write $R(x,y)$ for $R^W(x,y)=R^F(x,y)$ in the Theorem~\ref{thm:R(x,y)-on-Zd}.
\end{remark}

\begin{theorem}\label{thm:R(x,y)-on-Zd}
  Resistance distance on the integer lattice $(\bZd,\one)$ is given by
  \linenopax
  \begin{equation}\label{eqn:R(x,y)-on-Zd}
    R(x,y) = \frac1{(2\gp)^d} \int_{\bT^d} \frac{\sin^2((x-y)\cdot \frac t2)} {S(t)} \,dt,
  \end{equation}
  where $S(t) = 4 \sum_{k=1}^d \sin^2\left(\frac{t_k}2\right)$ as in \eqref{eqn:S(t)}.
  \begin{proof}
    We compute the resistance distance via $R(x,y) = v_x(x) + v_y(y) - v_x(y) - v_y(x)$. Using $e_x = e^{\ii x \cdot t}$, substitute in the terms from \eqref{eqn:vx-on-Zd-as-exponential} of Lemma~\ref{thm:vx-on-Zd}:
    \linenopax
    \begin{align}
      R(x,y)
      &= \frac1{(2\gp)^d} \int_{\bTd} \frac{\cj{e_x}(e_x-1) + \cj{e_y}(e_y-1) - \cj{e_x}(e_y-1) - \cj{e_y}(e_x-1)}{S(t)} \,dt 
      %&= \frac1{(2\gp)^d} \int_{\bTd} \frac{1-\cancel{\cj{e_x}} + 1 - \cancel{\cj{e_y}} - e_{y-x} + \cancel{\cj{e_x}} - \cj{e_{y-x}} + \cancel{\cj{e_y}}}{S(t)} \,dt \notag \\
      %&= \frac1{(2\gp)^d} \int_{\bTd} \frac{2-2\cos((x-y) \cdot t)}{S(t)} \,dt, 
      \label{eqn:vx-on-Zd-cos}
    \end{align}
    and the formula follows by the half-angle identity and other algebra.
  \end{proof}
\end{theorem}

\begin{remark}\label{rem:monopoles}
  The fact that $(\bZd,\one)$ is transient if and only if $d \geq 3$ was first discovered by Polya \cite{Polya21}, and so Theorem~\ref{thm:monopoles-on-Zd} is a result which is well-known in the literature (cf.~\cite[Thm.~5.11]{Soardi94} and \cite{DoSn84,Nash-Will59}, e.g.). We include this result here because the present context allows for a brief proof which offers some intuition for this startling dichotomy.
\end{remark}

\begin{theorem}\label{thm:monopoles-on-Zd}
  The network $(\bZd,\one)$ has a monopole
  \linenopax
  \begin{equation}\label{eqn:monopoles-on-Zd}
    w_o(x) = \frac1{(2\gp)^d} \int_{\bTd} \frac{\cos(x \cdot t)}{S(t)}\,dt
  \end{equation}
  if and only if $d \geq 3$, in which case the monopole at $o$ is unique.
  \begin{proof}
    As in the proof of Lemma~\ref{thm:vx-on-Zd}, we use the Fourier transform to solve $\Lap w_o = \gd_o$ by converting it into $S(t) \hat w_o(t) = 1$.
    This gives \eqref{eqn:monopoles-on-Zd}, in which the integral converges because $\frac{\cos t}{S(t)} \approx \frac1{S(t)} \in L^1(\bTd)$ iff $d \geq 3$. To see this, note that upon switching to spherical coordinates, $1/S(\gr) = O(\gr^{-2})$, as $\gr \to 0$, and one requires
    \linenopax
    \begin{equation}\label{eqn:exm:infinite-lattices:spherical-integral}
      \frac1{S(t)} \in L^1(\bTd) \q \iff \q \int_0^1 |\gr^{-2}| \gr^{d-1} \, dS_{d-1} < \iy,
    \end{equation}
    where $dS_{d-1}$ is the usual $(d-1)$-dimensional spherical measure. Of course, \eqref{eqn:exm:infinite-lattices:spherical-integral} holds precisely when $-2+d-1>-1$, i.e., when $d \geq 3$.It remains to check that $w_o \in \HE$. 
    Applying \eqref{eqn:spectral-measures-Phi(gx)} with $\gx=\gd_o$, we obtain
    \linenopax
    \begin{align}\label{eqn:energy-of-monopole-on-Zd}
      \|w_o\|_\energy &= \|\gF\gd_o\|_\energy = \| \Lap^{-1/2} \gd_o\|_2
      %= \int_{\bTd} S(t) \hat w_o(t)^2\,dt
      = \int_{\bTd} \frac1{S(t)} \,dt
      < \iy,
    \end{align}
    by \eqref{eqn:exm:infinite-lattices:spherical-integral} again (and $\hat \gd_o = \one$, as noted above). 
    To see uniqueness, suppose $w'$ were another monopole. Then $\Lap(w_o-w') = \gd_o - \gd_o = 0$ and $w_o-w'$ is harmonic. By Remark~\ref{rem:Harm-on-Zd}, the only finite-energy harmonic functions on $(\bZd,\one)$ are constant.
  \end{proof}
\end{theorem}

\begin{remark}
  Comparing \eqref{eqn:monopoles-on-Zd} to \eqref{eqn:vx-on-Zd} gives a heuristic as to why all networks support finite-energy dipoles, but not all support monopoles: the numerator in the integral for the monopole is $o(1)$ as $t \to 0$, while the corresponding numerator for the dipole is $o(t)$ as $t \to 0$.  
\end{remark}

\begin{cor}\label{thm:vx-in-L2(Zd)-for-big-d}
  For $(\bZd,\one)$, one has $v_x \in \ell^2(\bZd)$ if and only if $d \geq 3$.
  \begin{proof}
    %By computations similar to those in the proof of Theorem~\ref{thm:finite-resistance-to-infinity-in-Zd}, 
    One can see that in absolute values, the integrand $\left|(e^{\ii x\cdot t} - 1)/S(t)\right|$ of \eqref{eqn:vx-on-Zd-as-exponential} is in $L^2(\bTd)$ if and only if $d \geq 3$ (one only needs to check for $t \approx 0$, which is easy in spherical coordinates), in which case Parseval's theorem applies.
  \end{proof}
\end{cor}

Corollary~\ref{thm:vx-in-L2(Zd)-for-big-d} is a result comparable to \cite[Prop.~2]{CaW92}, but for a situation in which the isoperimetric inequality is not satisfied; see \eqref{eqn:strong-isoperimetric-inequality}.

\begin{cor}\label{thm:w-in-L2(Zd)-for-big-d}
  For $(\bZd,\one)$, the monopoles $w_x$ lie in $\ell^2(\bZd)$ if and only if $d \geq 5$.
  \begin{proof}
    The proof is the same as in Corollary~\ref{thm:vx-in-L2(Zd)-for-big-d}, except that the integrand is $1/S(t)$, which is in $L^2(\bTd)$ if and only if $d \geq 5$.
  \end{proof}
\end{cor}

%\subsection{Spectral resolution of \Lap on $(\bZd,\one)$}

We now apply some of our results from previous sections to obtain an explicit formula for the spectral resolution of $\Lap_2$ and \LapF, for the example $(\bZd,\one)$. 

\begin{lemma}\label{thm:R(x,y)-in-spectral-fourier}
  For $(\bZd,\one)$, if $\gx = \gd_x - \gd_y$ for any fixed $x,y \in G$, then the corresponding spectral measure of $\Lap_2$ is
  \linenopax
  \begin{align}\label{eqn:measure-conversion}
    \gl \left( \int_{S(t)=\gl} (1-\cos((x-y)\cdot t) \,dt \right) d\gm_\gx^{\ell^2} 
    = \frac1{(2\gp)^d} (dt) \comp S^{-1},
  \end{align}
  where $S(t) = 4 \sum_{k=1}^d \sin^2\left(\frac{t_k}2\right)$ as in \eqref{eqn:S(t)}.
  \begin{proof}
    This is a standard result from spectral theory, but we include it because explicit formulas are not commonplace in this subject. Since $\Lap_2$ corresponds to multiplication by \gl on the spectral side and multiplication by $S(t)$ on the Fourier side, it is easiest to see \eqref{eqn:measure-conversion} from  
    \linenopax
    \begin{align}\label{eqn:R(x,y)-conversion}
      \int_0^\iy \frac1\gl \|P_2(d\gl)(\gd_x - \gd_y)\|_2^2
      &= \frac1{(2\gp)^d} \int_{\bT^d} \frac{\sin^2((x-y)\cdot \frac t2)} {S(t)} \,dt,
    \end{align}
    which follows by comparing the two expressions for $\|\Lap_2^{-1/2}(\gd_x-\gd_y)\|_2^2$ found in Corollary~\ref{thm:R(x,y)} and in Theorem~\ref{thm:R(x,y)-on-Zd}.
%    \linenopax
%    \begin{align}\label{eqn:sim}
%      \|P_2(d\gl)(\gd_x - \gd_y)\|_2
%      \sim \frac{\sin^2((x-y)\cdot \frac t2)}{(2\gp)^{d/2}}.
%    \end{align}
%    More precisely, the Fourier transform of $\gd_x(n) \in \ell^2(\bZd,\one)$ is $\hat \gd_x(t) = e^{\ii x \cdot t}$, and for each Borel set $A$, the spectral projection $P_2(A)$ corresponds to multiplication by \charfn{A}, whence the multiplication operator $M_{S(t)}$ corresponds to \charfn{S^{-1}(A)}.
  \end{proof}
\end{lemma}

\subsection*{Acknowledgements}

The authors are grateful for stimulating comments, helpful advice, and valuable references from Dorin Dutkay, Daniel Lenz, Paul Muhly, Christian Remling, and others. We also thank the students and colleagues who have endured our talks on this material and raised fruitful questions.

{\small
\bibliographystyle{alpha}%model1-num-names}%alpha}
\bibliography{networks}
}

\end{document}